\documentclass[12pt,a4paper]{article}
\usepackage[T2A]{fontenc}
\usepackage{srcltx}\usepackage{graphicx,color}
\usepackage[all]{xy}
\usepackage{latexsym,amsfonts,amssymb,amsmath,longtable,
mathrsfs,theorem,cite,hyperref
}
\newcommand{\D}{\mathscr D}

\newcommand{\qed}{\hfill{$\Box$}}
\theoremstyle{plain}
{

\newtheorem{Lemma}{{\bfseries Lemma}}[section]
{\theorembodyfont{\itshape}

\newtheorem{Theo}{{\bfseries Theorem}}
\newtheorem{Prop}{{\bfseries Proposition}}
\newtheorem{Pro}{{\bfseries Proposition}}[section]
\newtheorem{Probl}{{\bfseries Problem}}
\newtheorem{Cor}{{\bfseries Corollary}}[Theo]}
\theorembodyfont{\rm}

\DeclareMathOperator{\Inn}{Inn} 
 
\DeclareMathOperator{\Aut}{Aut}

\DeclareMathOperator{\Syl}{Syl} \DeclareMathOperator{\Hall}{Hall}

\DeclareMathOperator{\m}{m} \DeclareMathOperator{\sm}{sm}

\setlength{\textwidth}{160mm}
\setlength{\textheight}{250mm} \headheight0mm \headsep0mm
\oddsidemargin0mm \topmargin -10mm

\title{\vspace{-1cm} \hfill{\normalsize 20D60}{
\fontfamily{cmr} \fontseries{bx} \selectfont \\ \vspace{1cm} 
Classification and properties\\ of the $\pi$-submaximal subgroups\\ in minimal nonsolvable groups}
\thanks{The first  author is supported by
a NNSF grant of China 
and Wu Wen-Tsun Key Laboratory of Mathematics of Chinese Academy of Sciences. The second  author is supported by Chinese Academy of Sciences President's International Fellowship Initiative, PIFI, (Grant \#~2016VMA078).}}
\date{}

\author{Wenbin Guo\\
{\small Department of Mathematics, University of Science and
Technology of China,}\\ {\small Hefei 230026, P. R. China}\\
{\small E-mail:
wbguo@ustc.edu.cn}\\ \\
Danila O. Revin\\
{\small 1. Department of Mathematics, University of Science and
Technology of China,}\\ {\small Hefei 230026, P. R. China}\\
{\small 2. Sobolev Institute of Mathematics SB RAS, and}\\
{\small 3.  Novosibirsk State University,}\\
{\small Novosibirsk 630090, Russia}\\
{\small E-mail: revin@math.nsc.ru}}

\begin{document}

\maketitle
\pagenumbering{arabic}
\begin{abstract}
Let $\pi$ be a set of primes. According to H.~Wielandt, a subgroup $H$ of a 
finite group $X$ is called a {\it $\pi$-submaximal subgroup} if there is a monomorphism $\phi:X\rightarrow Y$ into a finite group 
$Y$ such that $X^\phi$ is subnormal in $Y$ and $H^\phi=K\cap X^\phi$ for a $\pi$-maximal subgroup $K$ of $Y$. 
In his talk at the well-known conference on finite groups in Santa-Cruz (USA) in 1979, Wielandt posed a series of open questions and among them the 
following problem: to describe the $\pi$-submaximal subgroup of the minimal nonsolvable groups and to study properties of such subgroups: 
the pronormality, the intravariancy, the conjugacy in the automorphism group etc. In the article, for every set $\pi$ of primes, we obtain a description of the $\pi$-submaximal subgroup in minimal nonsolvable groups and investigate their properties, so we give a solution of Wielandt's problem.

\smallskip

\noindent{\bf Key words:} Minimal nonsolvable group, minimal simple group, $\pi$-maximal subgroup, $\pi$-submaximal subgroup, pronormal subgroup.
\end{abstract}


\section{Introduction}

In the paper, all groups are finite and $G$ always denotes a finite group. Moreover, $\pi$~denotes some given set of primes and $\pi'$ is the set of primes $p$ such that $p\notin\pi$. For a natural number $n$, let $\pi(n)$ be the set of prime divisors of $n$ and let $n_\pi$ be the {\it $\pi$-share} of~$n$, that is, the greatest divisor $m$ of $n$ such that $\pi(m)\subseteq\pi$. Clearly, $n=n_\pi n_{\pi'}$ and $(n_\pi, n_{\pi'})=1$. For a group $G$, let $\pi(G)=\pi(|G|)$. A  group $G$ is called a {\it $\pi$-group} if $\pi(G)\subseteq\pi$.

A subgroup $H$ of  $G$ is said to be {\it $\pi$-maximal} if $H$ is maximal with respect to inclusion in the set of 
$\pi$-subgroups  of~$G$. We denote by~$\m_\pi(G)$ the set of $\pi$-maximal subgroups of $G$.

A group $G$ is called a {\it $\D_\pi$-group} (or write $G\in \D_\pi$) if any two subgroups in $\m_\pi(G)$ are conjugate.  
It is well-known that $G\in\D_\pi$ if and only if the complete analog of the Sylow theorem holds for $\pi$-subgroups of $G$, that is,
\begin{itemize}
  \item[$(1)$] $G$ possesses a  {\it $\pi$-Hall} subgroup (that is, a subgroup of order~ $|G|_\pi$),
  \item[$(2)$] every two $\pi$-Hall subgroups are conjugate in~$G$,
  \item[$(3)$] every $\pi$-subgroup of $G$ is contained in some $\pi$-Hall subgroup.
\end{itemize}
In particular, in a $\D_\pi$-group, the $\pi$-maximal subgroups are exactly the $\pi$-Hall subgroups. The Sylow theorem states that $G\in\D_p$ for any group $G$ and any prime $p$. The well-known Hall--Chunikhin theorem \cite{Chun,Hall1,Hall2} states that a group $G$ is solvable if and only if $G\in\D_\pi$ for every set $\pi$ of primes. But in general, for a given set $\pi$, the class of $\D_\pi$-groups is wider than the class of solvable groups since every $\pi$-group is a $\D_\pi$-group.

We denote by $\Hall_\pi(G)$ the set of $\pi$-Hall subgroups of a group $G$. 
The $\D_\pi$-groups and the groups containing $\pi$-Hall subgroups are  well investigated (see a survey~\cite{Surv} and next 
works~\cite{GRV,FrattArg,GR,ExProHall}). In particular, it is proved that $G\in \D_\pi$ if and only if every composition factor of $G$ is a $\D_\pi$-groups (see \cite[Theorem~7.7]{Hall3'}, \cite[Theorem~6.6]{Surv}, \cite[Chapter~2, Theorem~6.15]{GuoBook1}) and the simple $\D_\pi$-groups are described in~\cite[Theorem~3]{R4} (see also~\cite[Theorem~6.9]{Surv}). These results are based on the description of Hall subgroups in the finite simple groups and the following nice properties of $\pi$-Hall subgroups: if $N$ is  normal and $H$ is  $\pi$-Hall subgroups in $G$, then $HN/N\in\Hall_\pi(G/N)$ and $H\cap N\in\Hall_\pi(N)$.

The Hall--Chunikhin theorem shows that, if $G$ is nonsolvable, then $G\notin \D_\pi$ 
for some $\pi$. It is clear that $\m_\pi(G)\ne\varnothing$ and $\Hall_\pi(G)\subseteq\m_\pi(G)$ for any group $G$ but, in general,  $\m_\pi(G)\nsubseteq\Hall_\pi(G)$ and it may happen that $\Hall_\pi(G)=\varnothing$.
 It is natural to try  to find a description for $\pi$-maximal subgroups similar to ones for $\pi$-Hall subgroups. A hardness is that the  $\pi$-maximal subgroups have no properties similar to the mentioned above properties of $\pi$-Hall subgroups.

For example \cite[(4.2)]{Wie3}, if $X$ and $Y$ are groups and $X\notin\D_\pi$ for some $\pi$, then, for any $\pi$-subgroup (not only for $\pi$-maximal) $K$ of $Y$, there is a $\pi$-maximal subgroup $H$ of the regular wreath product $X\wr Y$ such that the image of $H$ under the natural epimorphism $X\wr Y\rightarrow Y$ coincides with $K$.

In general, the intersection of a  $\pi$-maximal subgroup $H$ with a normal subgroup $N$ of  $G$ is not a $\pi$-maximal subgroup of~$N$.  For example, it is easy to show that a Sylow 2-subgroup $H$ of $G=PGL_2(7)$ is  $\{2,3\}$-maximal in $G$ but $H\cap N\notin\m_{\{2,3\}}(N)$ for $N=PSL_2(7)$.

But the situation with intersections of $\pi$-maximal and normal 
subgroups of a finite group is not so dramatic in comparing with the situation with images under homomorphisms since the following statement holds.

\begin{Prop}\label{WieHart} {\rm (Wielandt--Hartley Theorem)  }{\it Let $G$~be a finite group, let $N$~be a subnormal subgroup of~$G$ and  $H\in \m_\pi(G)$. Then  $H\cap N=1$ if and only if $N$ is a $\pi'$-group.}
\end{Prop}

For the case where $N$ normal in $G$, Wielandt's proof of this statement can be found in \cite[13.2]{Wie4}, and Hartly's one in \cite[Lemmas 2 and~3]{Hart}. The proof of the above proposition in the general case see \cite[Theorem~7]{Shem}.

In view of Proposition~\ref{WieHart}, it is natural to consider the following concept. According to H.~Wielandt, 
a subgroup $H$ of a group $X$ is called a {\it $\pi$-submaximal subgroup} if there is a 
monomorphism $\phi:X\rightarrow Y$ into a  group $Y$ such that $X^\phi$ is subnormal in $Y$ and $H^\phi=K\cap X^\phi$ for a $\pi$-maximal 
subgroup $K$ of $Y$. We denote by~$\sm_\pi(X)$ the set of $\pi$-maximal subgroups of~$X$.

Evidently, $\m_\pi(G)\subseteq\sm_\pi(G)$ for any group $G$. 
The inverse inclusion does not hold in general as one can see in the above example: any Sylow 2-subgroup of $PSL_2(7)$ is $\{2,3\}$-submaximal but is not   $\{2,3\}$-maximal.

In view of the Wielandt--Hartley theorem, $1\in \sm_\pi(G)$ if and only if $G$ is  a $\pi'$-group\footnote{Note that, in Wielandt's talk~\cite[(5.4)(a)]{Wie3}, one can find a more strong version of this statement: {\it if $H\in\sm_\pi(G)$ then $N_G(H)/H$ is a $\pi'$-group.} For the case where $G$ is simple see Lemma~\ref{WieHartSimple}.}. Moreover, in contrast with $\pi$-maximal subgroups,
the following property of $\pi$-submaximal subgroups holds:
\begin{description}
\item If $H\in  \sm_\pi(G)$ and $N$ is a (sub)normal subgroup of $G$, then $H\cap N\in\sm_\pi(N)$. \ \ \ $(*)$
\end{description}
This property shows that, in some sense, the behavior of  $\pi$-submaximal subgroups is similar to ones of $\pi$-Hall 
subgroups under taking of intersections with (sub)normal subgroups. Note that properties of $\pi$-submaximal subgroups are investigated in~\cite{Wie2,Wie4,GR1} (see also~\cite{Wie3,Hart1}).

By using the closeness of the class of $\D_\pi$-group under taking 
extensions \cite[theorem~7.7]{Hall3'}, (see also\cite[theorem~6.6]{Surv} and \cite[Chapter~2, Theorem~6.15]{GuoBook1}) 
one can show that $G\in\D_\pi$ if and only if every two $\pi$-submaximal subgroup of $G$ are conjugate. In particular,  
if $G\in\D_\pi$, then $\sm_\pi(G)=\m_\pi(G)=\Hall_\pi(G)$. 

In view of the Hall--Chunikhin theorem, it is natural to consider some ``critical'' situation where $G$ is non-solvable (and  
$G\notin \D_\pi$ for some $\pi$) but all subgroups of $G$ are solvable (and so they are $\D_\pi$-groups). In this situation, 
$G$ possesses more than one conjugacy
class of $\pi$-submaximal subgroups. 
In the paper, we consider the following problem which was possed by H.~Wielandt in his talk\footnote{This problem is one of ten open problems 
in this talk~\cite{Wie3}. Another problem \cite[Frage (i)]{Wie3} 
is the following conjecture: {\it a subgroup $A$ of a group $G$ is subnormal if $N_A(H\cap A)/(H\cap A)$ is $\pi'$-group for every set $\pi$ and 
every $H\in \m_\pi(G)$.} This is a converse statement to the strong version of the Wielandt--Hartley theorem (see previous remark) and it was proved by 
P.~Kleidman in his well-known work~\cite{Kleid}.} at the well-known conference on finite groups in Santa-Cruz  in 1979 \cite[Frage (g)]{Wie3}:

\begin{Probl}\label{WieProb}  {\rm (H.~Wielandt, 1979) To describe the $\pi$-submaximal subgroups of the minimal nonsolvable groups. To study properties of such subgroups: conjugacy classes, the pronormality, the intravariancy\footnote{In Wielandt's terminology, a subgroup $H$ of a group $G$ is said to be   {\it intravariant} if the conjugacy class of $H$ in $G$ is $\Aut(G)$-invariant, that is, for any $\alpha \in \Aut(G)$, there is $g\in G$ such that $H^\alpha=H^g$. Recall,  a subgroup $H$ of a group $G$ is said to be   {\it pronormal} if $H$ and $H^g$ are conjugate in $\langle H,H^g\rangle$ for any $g\in G$.}, the conjugacy in the automorphism group etc.}
\end{Probl}

A group $G$ is {\it minimal nonsolvable} if $G$ is nonsolvable but every proper subgroup of $G$ is solvable. It is well-known that $G$ is a minimal nonsolvable group if and only if $G/\Phi(G)$ is a  finite minimal simple group (that is, a non-abelian finite simple group $S$ such that $S$ is minimal nonsolvable), where $\Phi(G)$ is the Ftattini subgroup of~$G$ that is the intersection of the maximal subgroups of~$G$. In 1968, J.~Thompson \cite[Corollary~1]{Thomp} proved that $S$ is a minimal simple  group if and only if $S$ is isomorphic to a group in the following list $\cal T$ (we will call it as the Thompson list):

$(1)$ $L_2(2^p)$ where $p$~is a prime;

$(2)$ $L_2(3^p)$ where $p$~is an odd prime;

 $(3)$ $L_2(p)$ where $p$~is a prime such that $p>3$ and $p^2+1\equiv 0\pmod 5$;

  $(4)$ $Sz(2^p)$ where $p$~is an odd prime;

  $(5)$ $L_3(3)$.

  Actually, the properties emphasized by Wielandt in Problem~\ref{WieProb} play an important role in the study of $\pi$-maximal and, in particular, $\pi$-Hall subgroups. For a $\pi$-Hall subgroup $H$ of  $G$, the intravariancy means that $H$ can be lifted to a $\pi$-Hall subgroup $K$ of $\Aut(G)$ such that $K\cap \Inn(G)$ coincides with the image of $H$ in $\Inn(G)$ under the natural epimorphism $G\rightarrow\Inn(G)$ \cite[Proposition~4.8]{Surv}. One of important results concerning with $\pi$-Hall subgroups of the finite simple groups is the statement that if a finite simple group $S$ contains a $\pi$-Hall subgroup, then the number of the conjugacy classes of $\pi$-Hall subgroups in $S$ is a bounded $\pi$-number \cite[Theorem~1.1]{NumbCl}, \cite[Theorem~3.4]{Surv}. This statement has many consequences for arbitrary finite groups: criteria of the existence \cite{ExCrit} and the conjugacy \cite{ConjCrit}, closeness of the class of $\D_\pi$-groups under taking extensions, some analog of the Frattini argument 
for Hall subgroups~
   \cite{FrattArg}~etc. The strongest results on Hall subgroups turned out to be formulated in term of the pronormality \cite{ProSimple, ProHall, ExProHall, GR}: the pronormality of Hall subgroups in the finite simple groups, the existence of a pronormal $\pi$-subgroup in any group containing $\pi$-Hall subgroups etc.

In the paper, for every set $\pi$ of primes and for every  minimal nonsolvable group $G$, we solve Problem~\ref{WieProb}. More precisely, in the first, we reduce this problem to the case where $G$ is a  minimal simple group by proving the following statement:

\begin{Prop}\label{SubFactor} {\it Let $\pi$ be a set of primes.   Let $G$~be a finite group and let $N=F(G)$~be the Fitting subgroup of $G$ (that is, the greatest normal nilpotent subgroup of~$G$). Suppose $H\in\sm_\pi(G)$. Then the following statements hold:
\begin{itemize}
\item[$(1)$]   $\m_\pi(G/N)=\{KN/N\mid K\in \m_\pi(G)\}$.
  \item[$(2)$]  ${H\cap N}$ coincides with the $\pi$-Hall subgroup $O_\pi(N)$ of $N$.
  \item[$(3)$]   $HN/N\in\sm_\pi(G/N)$.
  \item[$(4)$]  
  $HN/N$ is pronormal in $G/N$ if and only if $H$ is pronormal in~$G$.
  \item[$(5)$]  $H$ is intravariant in~$G$ if and only if the conjugacy class of $HN/N$ in~$G/N$  is invariant under the image  ${\overline{\Aut(G)}}$ of the map $\Aut(G)\rightarrow\Aut(G/N)$ given by the rule $\phi\mapsto\bar\phi$ where $\bar\phi:Ng\mapsto Ng^\phi$ for $\phi\in \Aut(G)$.
\end{itemize}
}
\end{Prop}

If $G$ is a minimal nonsolvable group, then $F(G)=\Phi(G)$. Thus, for every $\pi$-submaximal subgroups $H$ of a minimal nonsolvable group $G$, the image of $H$ in the corresponding minimal simple group $\tilde{G}=G/\Phi(G)$ is a $\pi$-submaximal subgroup of $\tilde{G}$. In order to solve Problem~\ref{WieProb}, we only need to describe the $\pi$-submaximal subgroup of the minimal simple groups, that is, in the groups of the Thompson list. It is well-known 
\cite[Corollary~1.7.10]{GuoBook} that $\pi(G)=\pi(\tilde{G})$, and it is clear that, if $\pi\cap\pi(\tilde{G})=\varnothing$, then $\sm_\pi(G)=\m_\pi(G)=\{1\}$; if $\pi\cap\pi(\tilde{G})=\{p\}$, then $\sm_\pi(G)=\m_\pi(G)=\Syl_p(G)$; and if $\pi(\tilde{G})\subseteq\pi$, then $\sm_\pi(G)=\m_\pi(G)=\{G\}$.
A description of  the $\pi$-submaximal subgroup in the minimal simple groups for the remaining cases is given in the following theorem.

\begin{Theo}\label{Th1}{\it Let $\pi$ be a set of primes and $S$  a minimal simple group. Suppose that $|{\pi\cap\pi(S)}|>1$ and $\pi(S)\nsubseteq\pi$. Then  representatives of the conjugacy classes of $\pi$-sub\-maxi\-mal subgroups of~$S$, the information of their structure, $\pi$-maximality, pronormality, intravariancy, and the action of $\Aut(S)$ on the set of conjugacy classes of $\pi$-submaximal subgroups can be specify in the corresponding table~{\rm\ref{L_2(2),2_notin_pi}--\ref{L_3(3),2,3}} below.}
\end{Theo}

\begin{Cor}\label{Cor}{\it For every set $\pi$ of primes, the $\pi$-submaximal subgroups of any minimal nonsolvable group are pronormal.  }
\end{Cor}

Note that our results depend on Thompson's classification of the minimal simple groups and do not depend on the classification of finite simple  groups.

 \bigskip

 \noindent{\large\bf \thesubsection}\hspace{10pt} Notation in Tables~\ref{L_2(2),2_notin_pi}--\ref{L_3(3),2,3}
 
 \bigskip

\noindent According to \cite{Bray,Atlas, KL}, we use the following notation.

\begin{itemize}
\item[] $\varepsilon$  denotes either $+1$ or $-1$ and the sing of this number.

  \item[]$C_n$ denotes the cyclic group of order $n$.

  \item[]$E_q$  denotes the elementary abelian group of  order $q$ where $q$ is a power of a prime.

  \item[]$D_{2n}$  denotes the dihedral group of order $2n$, i.~e. \\$D_{2n}=\langle x,y\mid x^n=y^2=1, x^y=x^{-1}\rangle$. Note that $D_4\cong E_4$.

  \item[]$SD_{2^n}$  denotes the semi-dihedral group of order $2^n$, i.~e.\\$SD_{2^n}=\langle x,y\mid x^{2^{n-1}}=y^2=1, x^y=x^{2^{n-2}+1}\rangle$.

  \item[]$S_n$ means a symmetric group of degree $n$.

  \item[]$A_n$  denotes the alternating group of degree $n$.

  \item[]$GL_n(q)$  denotes the general linear group of degree $n$ over the field of order~$q$.

  \item[]$PGL_n(q)$  denotes the projective general linear group of degree $n$ over the field of order~$q$.

  \item[]$L_n(q)=PSL_n(q)$  denotes the projective special linear group of degree $n$ over the field of order~$q$.

  \item[]$r_+^{1+2n}$  denotes the extra special group of  order $r^{1+2n}$ and of exponent $r$ where $r$ is an odd prime.

  \item[]$A:B$ means a split extension of a group $A$ by a group $B$ ($A$ is normal).

  \item[]$A^n$  denotes the direct product of $n$ copies of $A$.

  \item[]$A^{m+n}$ means $A^m:A^{n}$.

  \item[] The conditions in the column ``Cond.'' are necessary and sufficient for the existence and the $\pi$-submaximality of corresponding $H$. If a cell in this column is empty, then it means that the  corresponding $\pi$-submaximal subgroup always exists.

 \item[] In the column ``$H$'' the structure of corresponding $H$ is given.

 \item[] The conditions in the column ``is not $\pi$-max. if'' are necessary and sufficient for corresponding $H$ to be not $\pi$-maximal in $S$. If either this column is skipped or a cell in this column is empty, then the corresponding subgroup is $\pi$-maximal.

  \item[] A number $n$ in the column ``NCC'' is equal to the number of conjugacy classes of $\pi$-submaximal subgroups of $S$ isomorphic to corresponding subgroup $H$ and, if $n>1$, then in the same column the action of $\Aut(S)$ on these classes is described.

  \item[] The symbol ``\checkmark '' in the column ``Pro.'' means that the corresponding subgroup $H$ is pronormal in $S$.

  \item[] The symbol ``\checkmark '' in the column ``Intra.'' means that the corresponding subgroup $H$ is intravariant in $S$. If a cell in this column is empty, then $H$ is not intravariant.

\end{itemize}

\subsection{\normalsize\rm The $\pi$-submaximal subgroups in $S=L_2(q)$, where $q=2^p$, $p$ is prime,\\ for $\pi$ such that $|\pi\cap\pi(S)|>1$ and $\pi(S)\nsubseteq\pi$ }

$$|S|=q(q-1)(q+1),\,\,\,\pi(S)=\{2\}\cup\pi(q-1)\cup\pi(q+1).$$


\begin{longtable}{|c|c|c|c|c|c|}\caption{{\it The $\pi$-submaximal subgroups of $S=L_2(q)$, where
$q=2^p$, $p$ is a prime. Case: $2\notin\pi$}\\
~\\
Notation: $\pi_\varepsilon=\pi\cap\pi(q-\varepsilon)$, $\varepsilon\in\{+,-\}$.}
\label{L_2(2),2_notin_pi}\\  \hline

&Cond. & $H$ & NCC  & Pro. & Intra. \\
\hline\hline

1& $\pi_+\ne \varnothing$ & $C_{(q-1)_\pi}$  & $1$   &  \checkmark & \checkmark       \\ \hline
2& $\pi_-\ne \varnothing$  & $C_{(q+1)_\pi}$  & $1$  &    \checkmark & \checkmark      \\ \hline

\end{longtable}
\begin{center}
In any cases $H$ is $\pi$-maximal.
\end{center}


\begin{longtable}{|c|c|c|c|c|c|}\caption{{\it The $\pi$-submaximal subgroups of $S=L_2(q)$, where
$q=2^p$, $p$ is a prime. Case:  $2\in\pi$}\\
~\\
Notation: $\pi_\varepsilon=\pi\cap\pi(q-\varepsilon)$, $\varepsilon\in\{+,-\}$
}
\label{L_2(2),2_in_pi}\\  \hline
&Cond. & $H$ & NCC &  Pro. & Intra. \\
\hline\hline

1&       & $E_q:C_{(q-1)_\pi}$  & $1$    &  \checkmark & \checkmark      \\ \hline
2& $\pi_+\ne \varnothing$  & $D_{2(q-1)_\pi}$  & $1$        &  \checkmark & \checkmark    \\
  \hline
3& $\pi_-\ne \varnothing$  & $D_{2(q+1)_\pi}$  & $1$     &  \checkmark & \checkmark     \\ \hline
 \end{longtable}
 \begin{center}
In any cases $H$ is $\pi$-maximal.
\end{center}

\subsection{\normalsize\rm The $\pi$-submaximal subgroups in $S=L_2(q)$, where $q=3^p$, $p$ is odd prime,\\ for $\pi$ such that  $|\pi\cap\pi(S)|>1$ and $\pi(S)\nsubseteq\pi$ }

$$|S|=\frac{1}{2}q(q-1)(q+1),\,\,\,\pi(S)=\{3\}\cup\pi(q-1)\cup\pi(q+1).$$


\begin{longtable}{|c|c|c|c|c|c|}\caption{{\it The $\pi$-submaximal subgroups of $S=L_2(q)$, where
$q=3^p$, $p$ is an odd prime. Case:  $2\notin\pi$}\\
~\\
Notation: $\pi_\varepsilon=\pi\cap\pi(q-\varepsilon)$, $\varepsilon\in\{+,-\}$.}
\label{L_2(3),2_notin_pi}\\  \hline

&Cond. & $H$ & NCC  & Pro. & Intra. \\
\hline\hline

 1& $3\in\pi$     & $E_q:C_{(q-1)_\pi}$  & $1$    &  \checkmark & \checkmark      \\ \hline
2& $3\notin\pi$ and $\pi_+\ne \varnothing$ & $C_{(q-1)_\pi}$  & $1$   &  \checkmark & \checkmark       \\ \hline
3& $\pi_-\ne \varnothing$  & $C_{(q+1)_\pi}$  & $1$  &    \checkmark & \checkmark      \\ \hline
\end{longtable}
\begin{center}
In any cases $H$ is $\pi$-maximal.
\end{center}


\begin{longtable}{|c|c|c|c|c|c|}\caption{{\it The $\pi$-submaximal subgroups of $S=L_2(q)$, where
$q=3^p$, $p$ is an odd prime. Case:  $2\in\pi$}\\
~\\
Notation: $\pi_\varepsilon=\pi\cap\pi(q-\varepsilon)$, $\varepsilon\in\{+,-\}$
}
\label{L_2(3),2_in_pi}\\  \hline
&Cond. & $H$ & NCC &  Pro. & Intra. \\
\hline\hline

1&    $3\in\pi$     & $E_q:C_{\frac{1}{2}(q-1)_\pi}$  & $1$    &  \checkmark & \checkmark      \\ \hline
2& $\pi_+\ne \{2\}$  & $D_{(q-1)_\pi}$  & $1$        &  \checkmark & \checkmark    \\
  \hline
 3&  & $D_{(q+1)_\pi}$  & $1$     &  \checkmark & \checkmark     \\ \hline
 4&  $3\in\pi$         & $A_4$            & $1$       &  \checkmark & \checkmark       \\
 \hline
 \end{longtable}
 \begin{center}
In any cases $H$ is $\pi$-maximal.
\end{center}

\subsection{\normalsize\rm The $\pi$-submaximal subgroups in $S=L_2(q)$, where $q$ is a prime, $q^2\equiv -1\pmod5$,\\ for $\pi$ such that  $|\pi\cap\pi(S)|>1$ and $\pi(S)\nsubseteq\pi$ }
$$|S|=\frac{1}{2}q(q-1)(q+1),\,\,\,\pi(S)=\{q\}\cup\pi(q-1)\cup\pi(q+1).$$


\begin{longtable}{|c|c|c|c|c|c|}\caption{{\it The $\pi$-submaximal subgroups of $S=L_2(q)$, where
$q>3$ is a  prime, $q^2\equiv -1\pmod5$. Case:  $2\notin\pi$}\\
~\\
Notation: $\pi_\varepsilon=\pi\cap\pi(q-\varepsilon)$, $\varepsilon\in\{+,-\}$.}
\label{L_2(q),2_notin_pi}\\  \hline

&Cond. & $H$ & NCC  & Pro. & Intra. \\
\hline\hline

1& $q\in\pi$        & $C_q:C_{(q-1)_\pi}$  & $1$    &  \checkmark & \checkmark    \\ \hline
2& $q\notin\pi$ and $\pi_+\ne \varnothing$ & $C_{(q-1)_\pi}$  & $1$   &  \checkmark & \checkmark       \\ \hline
3& $\pi_-\ne \varnothing$  & $C_{(q+1)_\pi}$  & $1$  &    \checkmark & \checkmark      \\ \hline

\end{longtable}
\begin{center}
In any cases $H$ is $\pi$-maximal.
\end{center}


\begin{longtable}{|c|c|c|c|c|c|c|}\caption{{\it The $\pi$-submaximal subgroups of $S=L_2(q)$, where
$q>3$ is a  prime, $q^2\equiv -1\pmod5$, in the case $2\in\pi$}\\
~\\
Notation: $\pi_\varepsilon=\pi\cap\pi(q-\varepsilon)$, $\varepsilon\in\{+,-\}$,\\ $\delta\in\Aut(S)\setminus\Inn(S)$, $|\delta|=2$,\\
$\Aut(S)=\langle\Inn(S),\delta\rangle\cong S:\langle\delta\rangle\cong PGL_2(q)$}
\label{L_2(q),2_in_pi}\\  \hline
&Cond. & $H$ & NCC & $H$  is not  $\pi$-max. if & Pro. & Intra. \\
\hline\hline

 1&$q\in\pi$        & $C_q:C_{\frac{1}{2}(q-1)_\pi}$  & $1$  & $\phantom{\displaystyle\frac{1}{2}}$  &  \checkmark & \checkmark      \\ \hline
 & either &            &     & either                          & &          \\
 & $\pi_+\ne\{2\}$, &            &      &   $\pi_+=\{2\}$, $3\in\pi$, and & & \\
 2& or  & $D_{(q-1)_\pi}$&$1$ &  $q\equiv
                                          41\pmod{48}{}^*$ & \checkmark  &   \checkmark   \\
& $3\notin\pi$,     &            &      &   or     &  &\\
&or &                  &      & $\pi_+=\{2,3\}$ and &    & \\
 &$q\equiv 1\pmod8$ &                  &      & $q\equiv7,31\pmod{72}{}^{**}$ &   & \\\hline
  & either &            &     & either                          & &          \\
 & $\pi_-\ne\{2\},$&         &      & $\pi_-=\{2\}$, $3\in\pi$, and  &  &     \\
3& or& $D_{(q+1)_\pi}$  & $1$  &    $q\equiv7
                                             \pmod{48}{}^*$    &  \checkmark & \checkmark     \\
&$3\notin\pi$, &          &      &     or                       &     &      \\
&or &                  &      & $\pi_-=\{2,3\}$ and &    & \\
 & $q\equiv-1\pmod8$ &                  &      & $q\equiv41,65\pmod{72}{}^{**}$ &   & \\ \hline
4& $3\in\pi$ and       & $A_4$            & $1$     &  
                                                                        &  \checkmark & \checkmark       \\
& $q\equiv\pm3\pmod8$&                 &            &             &             &            \\  \hline
5& $3\in\pi$ and        & $S_4$            & $2$ permuted     &   &  \checkmark &                \\
& $q\equiv\pm1\pmod8$&                 &  by $\delta$       &             &             &            \\ \hline
\end{longtable}
\begin{center}
${}^{*}$  $H\cong D_8$ and $H$ is contained in a $\pi$-maximal subgroup isomorphic to $S_4$;\\
${}^{**}$  $H\cong D_6$ and $H$ is contained in a $\pi$-maximal subgroup isomorphic to $S_4$.

\end{center}

\subsection{\normalsize\rm The submaximal $\pi$-subgroups in $S=Sz(q)$, where $q=2^p$, $p$ is odd prime,\\ for $\pi$ such that  $|\pi\cap\pi(S)|>1$ and $\pi(S)\nsubseteq\pi$ }

$$|S|=q^2(q-1)(q^2+1)=q^2(q-1)(q-r+1)(q+r+1),\text{ where }r=\sqrt{(2q)}=2^{(p+1)/2},$$

$$\pi(S)=\{2\}\cup\pi(q-1)\cup\pi(q- r+1)\cup\pi(q+ r+1),$$


\begin{longtable}{|c|c|c|c|c|c|}\caption{{\it The $\pi$-submaximal subgroups of $S=Sz(q)$, where
$q=2^p$, $p$ is an odd prime. Case:  $2\notin\pi$}\\
~\\
Notation: $r=\sqrt{(2q)}=2^{(p+1)/2}$\\ $\pi_0=\pi\cap\pi(q-1)$, $\pi_\varepsilon=\pi\cap\pi(q-\varepsilon r+1)$, $\varepsilon\in\{+,-\}$.}
\label{Sz,2_notin_pi}\\  \hline

&Cond. & $H$ & NCC  & Pro. & Intra. \\
\hline\hline

1&$\pi_0\ne \varnothing$ & $C_{(q- 1)_\pi}$  & $1$   &  \checkmark & \checkmark       \\ \hline
2& $\pi_+\ne \varnothing$ & $C_{(q- r+1)_\pi}$  & $1$   &  \checkmark & \checkmark       \\ \hline
3& $\pi_-\ne \varnothing$  & $C_{(q- r+1)_\pi}$  & $1$  &    \checkmark & \checkmark      \\ \hline

\end{longtable}
\begin{center}
In any cases $H$ is $\pi$-maximal.
\end{center}


\begin{longtable}{|c|c|c|c|c|c|}\caption{{\it The $\pi$-submaximal subgroups of $S=Sz(q)$, where
$q=2^p$, $p$ is an odd prime. Case:  $2\in\pi$}\\
~\\
Notation: $r=\sqrt{(2q)}=2^{(p+1)/2}$\\ $\pi_0=\pi\cap\pi(q-1)$, $\pi_\varepsilon=\pi\cap\pi(q-\varepsilon r+1)$, $\varepsilon\in\{+,-\}$.}

\label{Sz,2_in_pi}\\  \hline
&Cond. & $H$ & NCC &  Pro. & Intra. \\
\hline\hline

1&       & $E_q^{1+1}:C_{(q-1)_\pi}$  & $1$    &  \checkmark & \checkmark      \\ \hline
2&  $\pi_0\ne \varnothing$  & $C_{(q-1)_\pi}:C_4$  & $1$        &  \checkmark & \checkmark    \\
  \hline
3& $\pi_+\ne \varnothing$  & $D_{2(q- r+1)_\pi}$  & $1$        &  \checkmark & \checkmark    \\
  \hline
4& $\pi_-\ne \varnothing$  & $D_{2(q+ r+1)_\pi}$  & $1$     &  \checkmark & \checkmark     \\ \hline
 \end{longtable}
 \begin{center}
In any cases $H$ is $\pi$-maximal.
\end{center}

\subsection{\normalsize\rm The $\pi$-submaximal subgroups of $S=L_3(3)$,\\ for $\pi$ such that  $|\pi\cap\pi(S)|>1$ and $\pi(S)\nsubseteq\pi$ }

$$|S|=2^4\cdot3^3\cdot13,\,\,\, \pi(S)=\{2,3,13\},$$


\begin{longtable}{|c|c|c|c|c|}\caption{{\it The $\pi$-submaximal subgroups of $S=L_3(3)$.\\  Case:  $\pi\cap\pi(S)=\{3,13\}$}}

\label{L_3(3),3,13}
\\

\hline &$H$ & NCC  & Pro. & Intra. \\
\hline\hline

1& $C_{13}:C_{3}$  & $1$    &  \checkmark & \checkmark    \\ \hline
2& $3_+^{2+1}$  & $1$   &  \checkmark & \checkmark       \\ \hline

\end{longtable}
\begin{center}
In any cases $H$ is $\pi$-maximal.
\end{center}

\begin{longtable}{|c|c|c|c|c|}\caption{{\it  The $\pi$-submaximal subgroups of $S=L_3(3)$.\\ Case:  $\pi\cap\pi(S)=\{2,13\}$}}

\label{L_3(3),2,13}
\\

\hline &$H$ & NCC  & Pro. & Intra. \\
\hline\hline

1& $C_{13}$  & $1$    &  \checkmark & \checkmark    \\ \hline
2& ${SD}_{16}$  & $1$   &  \checkmark & \checkmark       \\ \hline

\end{longtable}
\begin{center}
In any cases $H$ is $\pi$-maximal.
\end{center}


\begin{longtable}{|c|c|c|c|c|c|}\caption{{\it  The $\pi$-submaximal subgroups of $S=L_3(3)$.\\ Case:  $\pi\cap\pi(S)=\{2,3\}$}\\
$\gamma\in \Aut(S)\setminus\Inn(S)$, $|\gamma|=2$,\\
$\Aut(S)=\langle\Inn(S),\gamma\rangle$}

\label{L_3(3),2,3}
\\  \hline
 &$H$ & NCC & $H$  is not  $\pi$-max. if & Pro. & Intra. \\
\hline\hline

 1&$E_{3^2}:GL_2(3)$  & $2$  permuted  by $\gamma$ &   &  \checkmark &   \\ \hline
 2&$3_+^{1+2}:C_{2}^2$   & $1$  &   always; $H\leq E_{3^2}:GL_2(3)$ &  \checkmark &    \checkmark \\ \hline
 3&$GL_2(3)$            & $1$     & always; $H\leq E_{3^2}:GL_2(3)$  &  \checkmark &  \checkmark  \\ \hline
 4&$S_4$            &  $1$    &   &  \checkmark &   \checkmark               \\\hline
\end{longtable}

\section{Preliminaries}

We write $M\lessdot G$ if $M$ is a maximal subgroup of $G$, that is, $M<G$ and $M
\leq H
\leq G$ implies that either $H=M$ or $H=G$. Moreover, we write $H\unlhd G$ and $H{\unlhd\unlhd}G$ if $H$ is a normal or subnormal subgroup of $G$, respectively.

\begin{Lemma}\label{MaxSubgrSimple} {\it
Let  $S$~be a minimal simple group. Then  representatives of the conjugacy classes of maximal subgroups of~$S$, the information on their structure, conjugacy classes,  and the action of $\Aut(S)$ on the set of the conjugacy classes of maximal subgroups can be specify in the corresponding table~{\rm\ref{MaxL_2(2)}--\ref{MaxL_3(3)}}.
}
\end{Lemma}

\noindent{\sc Proof}. See {\rm \cite[theorem~2.1.1, tables 8.1--8.4 and~8.16]{Bray}, \cite[theorem~II.8.27]{Hupp},  \cite{Atlas}, \cite[theorem~9]{SuzCl}}.\qed\medskip

\begin{longtable}{|c|c|}\caption{{\it Maximal subgroups of  $S=L_2(q)$\\ where $q=2^p$, $p$ is prime}
\label{MaxL_2(2)}}\\
 \hline

 $M$ & NCC  \\
\hline\hline

        $E_q:C_{(q-1)}$  & $1$    \\ \hline
  $D_{2(q-1)}$  & $1$           \\
  \hline
   $D_{2(q+1)}$  & $1$       \\ \hline
 \end{longtable}
 \begin{center}

\end{center}

\begin{longtable}{|c|c|}\caption{{\it Maximal subgroups of  $S=L_2(q)$\\ where
$q=3^p$, $p$ is an odd prime}}
\label{MaxL_2(3)}\\  \hline

 $M$ & NCC \\
\hline\hline

$C_q:C_{\frac{1}{2}(q-1)}$  & $1$     \\ \hline
 $D_{q-1}$  & $1$     \\
  \hline
$D_{q+1}$  & $1$   \\ \hline
 $A_4$            & $1$      \\
 \hline
 \end{longtable}

\begin{longtable}{|c|c|c|}\caption{{\it Maximal subgroups  of $S=L_2(q)$\\ where
$q>3$ is  prime, $q^2\equiv -1\pmod5$}}
\label{MaxL_2(q)}\\  \hline

 $M$ & NCC & Conditions\\
\hline\hline

 $C_q:C_{\frac{1}{2}(q-1)}$  & $1$  & $\phantom{\displaystyle\frac{1}{2}}$       \\ \hline
 $D_{(q-1)}$  & $1$   & $q\ne 7$    \\
 \hline
$D_{(q+1)}$  & $1$  &  $q\ne 7$   \\ \hline
$A_4$            & $1$     &  $q\equiv\pm3\pmod8$      \\ \hline

$S_4$            & $2$ permuted     &  $q\equiv\pm1\pmod8$              \\
                &  by $\delta$       &                   \\ \hline

\end{longtable}
\begin{center}

Here $\delta\in\Aut(S)\setminus\Inn(S)$, $|\delta|=2$,\\
$\Aut(S)=\langle\Inn(S),\delta\rangle\cong S:\langle\delta\rangle\cong PGL_2(q)$
\end{center}

\begin{longtable}{|c|c|}\caption{{\it Maximal subgroups of $S=Sz(q)$\\ where
$q=2^p$, $p$ is an odd prime}}
\label{MaxSz}\\  \hline

 $M$ & NCC  \\
\hline\hline

      $E_q^{1+1}:C_{(q-1)}$  & $1$       \\ \hline
 $C_{(q-1)}:C_4$  & $1$        \\
  \hline
 $D_{2(q- r+1)}$  & $1$     \\
  \hline
 $D_{2(q+ r+1)}$  & $1$         \\ \hline
 \end{longtable}


\begin{longtable}{|c|c|}\caption{{\it Maximal subgroups  of $S=L_3(3)$}\\
}
\label{MaxL_3(3)}
\\  \hline
 $M$ & NCC  \\
\hline\hline

 $E_{3^2}:GL_2(3)$  & $2$  permuted  by $\gamma$  \\ \hline
  $C_{13}:C_{3}$  & $1$   \\   \hline
 $S_4$            &  $1$               \\ \hline
\end{longtable}

\begin{center}
Here $\gamma\in \Aut(S)\setminus\Inn(S)$, $|\gamma|=2$,\\
$\Aut(S)=\langle\Inn(S),\gamma\rangle$
\end{center}

\begin{Lemma}\label{MaxSubgrAut} {\it
Let  $S$~be a minimal simple group, $M\lessdot G $ for some $G$ such that $\Inn(S)\leq G\leq \Aut(S)$, and $\Inn(S)\nleq M$. Then $G=MS$ and one of the following statement holds.
\begin{itemize}
  \item[$(1)$]  $M\cap \Inn(S)\lessdot \Inn(S)$;
  \item[$(2)$]  $S\cong L_2(7)$, $G=\Aut(S)$, $M\cong D_{12}$ and $M\cap \Inn(S)\cong D_6$;
  \item[$(3)$]  $S\cong L_2(7)$, $G=\Aut(S)$, $M\cong D_{16}$ and $M\cap \Inn(S)\cong D_8$;
  \item[$(4)$]  $S\cong L_3(3)$, $G=\Aut(S)$, $M\cong GL_2(3):C_2$ and $M\cap \Inn(S)\cong GL_2(3)$;
  \item[$(5)$]  $S\cong L_3(3)$, $G=\Aut(S)$, $M\cong 3_+^{1+2}:D_8$ and $M\cap \Inn(S)\cong 3_+^{1+2}:C_2^2$.
\end{itemize}
}
\end{Lemma}
\noindent{\sc Proof}. See {\rm \cite[theorem~2.1.1, tables 8.1--8.4 and~8.16]{Bray} and  \cite{Atlas}}.\qed\medskip

\begin{Lemma}\label{WieSubnorm} {\rm (see \cite[Statements 7 and~9]{Wie0}).}
Let $H_1,\dots,H_n$~be subnormal subgroups of~$G$. Then $K=\langle H_1,\dots,H_n\rangle$~is  subnormal in $G$ and every composition factor of $K$ is isomorphic to a composition factor of one of $H_1,\dots,H_n$.
\end{Lemma}

\begin{Lemma}\label{HallSubgroup}{\rm (see \cite[Lemma~1]{Hall} and \cite[Lemma 1.7.5]{GuoBook}).}
Let $A$ be a a normal subgroup and $H$ be a $\pi$-Hall subgroup of  $G$. Then
${H\cap A\in \Hall_\pi(A)}$ and  $HA/A\in \Hall_\pi(G/A).$
\end{Lemma}


\begin{Lemma}\label{NilpHall}{\rm (see \cite{Wie} and \cite[theorem~ 1.10.1]{GuoBook}).}
If $G$ has a nilpotent $\pi$-Hall subgroup, then $G\in\D_\pi$.
\end{Lemma}


Recall, a group $G$ is said to be $\pi$-separable if $G$ has a (sub)normal series
$$
G=G_0\unrhd G_1\unrhd\dots\unrhd G_n=1
$$
such that every section $G_{i-1}/G_i$ is either a $\pi$-group or a $\pi'$-group.

\begin{Lemma}\label{Separ}{\it If  $G$ is $\pi$-separable, then   $G\in \D_\pi$ and  $\sm_\pi(G)=\m_\pi(G)=\Hall_\pi(G)$.}
\end{Lemma}

\noindent{\sc Proof}. The statement $G\in\D_\pi$ is proved in \cite[Chapter~V, theorem~3.7]{Suz}. Thus, in view of $$\Hall_\pi(G)\subseteq\m_\pi(G)\subseteq\sm_\pi(G),$$ it is sufficient to prove  $\sm_\pi(G)\subseteq\Hall_\pi(G)$. Assume that $H\in\sm_\pi(G)$. Without loss of generality, we may assume that  $G{\unlhd\unlhd} X$ and $H=K\cap G$ for some $K\in\m_\pi(X)$. Let $Y=\langle G^X\rangle$ be the normal closure of $G$ in~$X$. Then $Y$ is $\pi$-separable by Lemma~\ref{WieSubnorm}. Moreover, $KY$ is $\pi$-separable,  so $KY\in\D_\pi$. Hence $K\in\m_\pi(HY)=\Hall_\pi(HY)$. Clearly, $G{\unlhd\unlhd}Y$. Hence $H=K\cap G\in\Hall_\pi(G)$ by Lemma~\ref{HallSubgroup}.
\qed\medskip

\begin{Lemma}\label{WieHartNorm}{\rm (see \cite[Lemma~2]{Hart}).} {\it Let  $A$~be a normal subgroup of~$G$  let $H\in \m_\pi(G)$. Then  $N_A(H\cap A)/(H\cap A)$ is a $\pi'$-group.}
\end{Lemma}

The following lemma was stated without a proof in \cite[5.3]{Wie3}. Here we give a proof for it.

\begin{Lemma}\label{SubmaxXSimple}\label{SubmaxSimple} {\it Let $\pi$ be a set of primes. Then, for a subgroup $H$ of a non-abelian simple group~$S$, the following conditions are equivalent.
\begin{itemize}
  \item[$(1)$]  $H\in \sm_\pi(S)$.
  \item[$(2)$] There exists a group $G$ such that $S$ is the socle of $G$,  $G/S$ is a $\pi$-group, and  $H=S\cap K$ for some $K\in\m_\pi(G)$.
  \item[$(3)$]  $H=S\cap K$ for some $K\in\m_\pi(\Aut(S))$ where $S$ is identified with $\Inn(S)$.
\end{itemize}}
\end{Lemma}

\noindent{\sc Proof}.

$(3)\Rightarrow(2):$ It is sufficient to take $G=SK$ where $K$ as in~$(3).$

$(2)\Rightarrow(1):$ It follows from the definition.

$(1)\Rightarrow(3)$. Take a group $X$ of the smallest order among all groups $G$ such that $S$ can be embedded into $G$ as a subnormal subgroup and $H=S\cap K$ for some $K\in\m_\pi(G)$. Let $A=\langle S^X\rangle$. Then the simplicity of $S$ implies that $A$~is a miniml normal subgroup in~$X$. Since  $X$ has the smallest possible order, we obtain $X=KA$.

Now we show that $K$ normalizes $S$, and as a consequence,  $A=\langle S^X\rangle=\langle S^K\rangle=S$. Clearly, $H=K\cap S\leq N_K(S)$. Hence $H\leq S\cap N_K(S)\leq S\cap K=H$ and $H=S\cap N_K(S)$. We claim  that $N=N_K(S)$~is a $\pi$-maximal subgroup of $X_0=SN$, from which we can obtain $K=N_K(S)$ and so $X=X_0\leq N_X(S)$. Assume that $N\leq U\leq X_0$ for a $\pi$-group~$U$. Since $SN\leq SU\leq X_0=SN$, we have $SU=SN$. Note that
$$U/(U\cap S)\cong US/S=NS/S\cong N/(N\cap S).$$ In order to prove that $N=U$, we only need to show that $U\cap S=N\cap S$. Let $U_0=U\cap S$ and let $g_1,\dots,g_m$~be a right transversal of $N$ in~$K$. Then the subgroups $S_i=S^{g_i}$, $i=1,\dots,m$, are pairwise different and $$A=\langle S_1,\dots, S_m\rangle\cong S_1\times\dots\times S_m.$$ Put $V=\langle U_0^{g_1},\dots, U_0^{g_m}\rangle.$ Then $K$ normalizes~$V$. Indeed, let $x\in K$. Since $K$ acts by the right multiplication on  $\{Ng_1,\dots,Ng_m\}$, there are a permutation~$\sigma$ of $\{1,\dots,m\}$ and some elements $t_1,\dots,t_m\in N$ such that $$g_ix=t_ig_{i\sigma},\,\,\,i=1,\dots,m.$$ Since $t_i\in N\leq U$ for all $i$, we have that $t_i$ normalizes $U_0=U\cap S$. Hence
\begin{multline*}
  V^x=
    \langle U_0^{g_1x},\dots, U_0^{g_mx}\rangle=
  \langle U_0^{t_1g_{1\sigma}},\dots, U_0^{t_mg_{m\sigma}}\rangle=\\
  \langle U_0^{g_{1\sigma}},\dots, U_0^{g_{m\sigma}}\rangle=
  \langle U_0^{g_{1}},\dots, U_0^{g_{m}}\rangle=V.
\end{multline*}
It follows that  $K$ normalizes $V$. Since $K\in\m_\pi(X)$, we have that  ${V\leq K}$. Consequently, $$U\cap S =U_0\leq V\cap S\leq K\cap S=H=N\cap S\leq U\cap S.$$ Thus, $U\cap S=N\cap S$ and  $U=N$.

 Now we have that $X=KS$, so $C_K(S)\unlhd X$ and $C_K(S)$ is a  $\pi$-group. Let $$\overline{\phantom{x}}:X\rightarrow X/C_K(S)$$ be the natural epimorphism. Note that its restriction to $S$ is an embedding of $S$ into the almost simple group  $\overline{X}$ with the socle $\overline{S}\cong S$.  It is easy to see that  $\overline{K}\in\m_\pi(\overline{X})$ and $\overline{H}=\overline{K}\cap \overline{S}$. By the choice of  $X$, we obtain $C_K(S)=1$,  so $X$~is almost simple.

 Thus, we can consider that $X\leq \Aut(S)$. Let $M\in\m_\pi(\Aut(S))$  such that
 $K\leq M$. Then $M\cap S\leq S\leq X$ and $K$ normalizes $M\cap S$. This implies that $M\cap S\leq K$ and $M\cap S=K\cap S =H$.
\qed\medskip

\begin{Lemma}\label{WieHartSimple} {\it Let $S$~be a finite simple group and  let $H\in \sm_\pi(S)$. Then  $N_S(H)/H$ is a $\pi'$-group.}
\end{Lemma}

\noindent{\sc Proof}. We identify $S$ with $\Inn(S)$. By Lemma~\ref{SubmaxSimple}, $H=K\cap S$, where $K\in \m_\pi(\Aut(S))$. Since $S\unlhd\Aut(S)$, Lemma~\ref{WieHartNorm} implies that $N_S(H)/H$ is  a $\pi'$-group. \qed\medskip

\begin{Lemma}\label{DpiInL2q} {\it Let $S=L_2(q)$ where $q$ is a power of some prime~$p$. Suppose that $\pi$ is a set of primes such that $2\notin\pi$ while $p\in\pi$. Then the following statements are equivalent.
\begin{itemize}
  \item[$(1)$] $S\in\D_\pi$.
  \item[$(2)$] $\Aut(S)\in\D_\pi$.
  \item[$(3)$] $\pi\cap\pi(S)\subseteq\{p\}\cup\pi(q-1)$.
\end{itemize}

}
\end{Lemma}

\noindent{\sc Proof}. See \cite[theorems A, 2.5, and,~3.3]{DpiCl}. \qed\medskip

\section{
Proof of Proposition~\ref{SubFactor} }
In this section, 
$N=F(G)$~is the Fitting subgroup of $G$ 
and $H\in\sm_\pi(G)$.

\begin{itemize}{\it
  \item[$(1)$]   $\m_\pi(G/N)=\{KN/N\mid K\in \m_\pi(G)\}$.}
\end{itemize}

   \noindent{\sc Proof}. Take $K\in\m_\pi(G)$. Suppose, $KN/N\leq L/N$ where $L/N$~is a $\pi$-group. Then ${L\in\D_\pi}$ by Lemma~\ref{Separ}. Since  $K\in\m_\pi(L)$, we have $K\in\Hall_\pi(L)$ and $KN/N\in\Hall_\pi(L/N)=\{L/N\}$ by Lemma~\ref{HallSubgroup}. Hence,  $KN/N\in\m_\pi(G/N)$ and $$\{KN/N\mid K\in \m_\pi(G)\}\subseteq\m_\pi(G/N).$$

   Conversely, assume that $L/N\in\m_\pi(G/N)$. By  Lemma~\ref{Separ}, $L\in\D_\pi$ and $L/N=KN/N$ for a fixed  $K\in\Hall_\pi(L)$. We show that $K\in\m_\pi(G)$. In fact, assume that $K\leq V$ for a $\pi$-subgroup $V$ of~$G$. Then $L\leq VN$ and $L=VN$ since $L/N\in\m_\pi(G/N)$. In view of  $K\in\Hall_\pi(L)$ and $L\in\D_\pi$, we have $K=V$ and $K\in\m_\pi(G)$.
 \qed\medskip


 \begin{itemize}{\it
  \item[$(2)$]  ${H\cap N}$ coincides with the $\pi$-Hall subgroup $O_\pi(N)$ of $N$.}
  \end{itemize}

 \noindent{\sc Proof}. Since $N$ is nilpotent, $N$ has a unique $\pi$-Hall subgroup. Hence we only need to show that ${{H\cap N}\in \Hall_\pi(N)}$. By  property $(*)$ of $\pi$-submaximal subgroups, $H\cap N\in\sm_\pi(N)$ and $\sm_\pi(N)=\Hall_\pi(N)$ by Lemma~\ref{Separ}.
  \qed\medskip

 \begin{itemize}{\it
  \item[$(3)$]   $HN/N\in\sm_\pi(G/N)$.}
  \end{itemize}

 \noindent{\sc Proof}. Denote by $\phi$ the natural epimorphism $G\rightarrow G/N$.
We need to show that  $H^\phi\in\sm_\pi(G^\phi)$. One can consider that there exists a group $X$ such that $G\,{\unlhd\unlhd}\, X$ and $H=G\cap K$ for some $K\in\m_\pi(X)$. Denote by $Y$ the normal closure $\langle N^X\rangle$ of $N$ in~$X$. Since $N\unlhd G\, {\unlhd\unlhd}\, X$, the Fitting theorem \cite[Chapter~A, theorem~8.8]{DH} implies that $Y$ is nilpotent. Consequently, $G\cap Y$ is also nilpotent. Since $G\cap Y\unlhd G$, we have
$$N\leq G\cap Y\leq F(G)=N.$$
Hence $G\cap Y=N$.
Consider the restriction $\tau:G\rightarrow X/Y$ to $G$ of the natural epimorphism $X\rightarrow X/Y$. Then the kernel of $\tau$ coincides with $N$. By the homomorphism theorem, there exists an injective homomorphism
$\psi: G^\phi=G/N \rightarrow X/Y$ such that the following diagram is commutative:
$$\xymatrix{
  G \ar[d]_{\phi} \ar[r]^{\tau} &     X/Y   \\
 G^\phi \ar[ur]_{\psi}                     }$$
Then  $$G^{\phi\psi}=G^\tau={GY/Y}\,{\unlhd\unlhd}\, {X/Y}.$$
Moreover, we have $$H^{\phi\psi}=H^\tau={HY/Y}=(G\cap K)Y/Y=(GY/Y)\cap (KY/Y)=G^{\phi\psi}\cap (KY/Y),$$  where $KY/Y\in\m_\pi(X/Y)$  in view of~(1). Thus $H^\phi\in\sm_\pi(G^\phi)$ by the definition of a $\pi$-sub\-maxi\-mal subgroup.
 \qed\medskip

 In order to prove next statements (4) and (5), we need the following lemma.

 \begin{Lemma}\label{SubmaxSepar}{\it In the above notation, $HN$ is $\pi$-separable and $H\in\Hall_\pi(HN)$.}
\end{Lemma}

\noindent{\sc Proof}. The group $HN$ has the subnormal series
   $$
  HN\unrhd N\unrhd H\cap N\unrhd 1
   $$
   such that  $HN/N$ is a $\pi$-group, $N/(H\cap N)$ is a $\pi'$-group by~(2) and $H\cap N$ is a $\pi$-group. Thus $HN$ is $\pi$-separable. Moreover, $H\in\Hall_\pi(HN)$ in view of
   $$|HN|_\pi=|HN/N||H\cap N|=|H/(H\cap N)||H\cap N|=|H|.$$
\qed\medskip

\begin{itemize}{\it
  \item[$(4)$] $HN/N$ is pronormal in $G/N$ if and only if $H$ is pronormal in~$G$.}
  \end{itemize}

   \noindent{\sc Proof}. If $H$ is pronormal in~$G$, then the pronormality of $HN/N$  in $G/N$ is evident.

   Conversely, assume that $HN/N$ is pronormal in $G/N$. Let $g\in G$. We need to prove that $H$ and $H^g$ are conjugate in $\langle H, H^g\rangle$.

   Firstly, consider the case when $g\in N_G(HN)$. Then $H^g\leq HN$ and $H^g\in\Hall_\pi(HN)$ by Lemma~\ref{SubmaxSepar}. Clearly, the subgroup $\langle H,H^g\rangle$ of the $\pi$-separable group $HN$ is also $\pi$-separable.  Lemma~\ref{Separ} implies that the $\pi$-Hall subgroups $H$ and $H^g$ of $\langle H,H^g\rangle$ are conjugate in $\langle H,H^g\rangle$.

   Now consider the general case for $g\in G$. Since $HN/N$ is pronormal in $G/N$, there is $y\in  \langle H,H^g\rangle$ such that $(HN)^y=(HN)^g$. Hence $gy^{-1}\in N_G(HN)$ and, in view of above, there exist some $z\in  \langle H,H^{gy^{-1}}\rangle\leq \langle H,H^g\rangle$ such that $H^z=H^{gy^{-1}}$. Hence $H$ and $H^g$ are conjugate by  $x=zy\in  \langle H,H^g\rangle$.
 \qed\medskip

   \begin{itemize}{\it
  \item[$(5)$]  $H$ is intravariant in~$G$ if and only if the conjugacy class of $HN/N$ in~$G/N$  is invariant under ${\overline{\Aut(G)}}$.}
  \end{itemize}

   \noindent{\sc Proof}. Recall that  ${\overline{\Aut(G)}}$ is the image in $\Aut(G/N)$ of $\Aut(G)$ under the map $\phi\mapsto\bar\phi$ where $\bar\phi:Ng\mapsto Ng^\phi$ for $\phi\in\Aut(G)$.

   If $H$ is intravariant in $G$, then for every $\phi\in\Aut(G)$, there is $g\in G$ such that  $H^\phi=H^g$. Hence $$(HN/N)^{\bar\phi}=\{(Nh)^{\bar\phi}\mid h\in H\}=\{Nh^{\phi}\mid h\in H\}=H^\phi N/N=H^gN/N$$
   and so the conjugacy class of $HN/N$ in~$G/N$  is invariant under ${\overline{\Aut(G)}}$.

   Conversely, assume that the conjugacy class of $HN/N$ in $G/N$ is invariant under ${\overline{\Aut(G)}}$. We need to show that for every $\phi\in \Aut(G)$ there exists $g\in G$ such that $H^\phi=H^g$. It is clear that $H^\phi\in\sm_\pi(G)$.
   By the hypothesis, there is $x\in G$ such that $H^\phi N=H^x N$, and in view of Lemma~\ref{SubmaxSepar}, $H^\phi ,H^x \in \Hall_\pi(H^\phi N)$ and $H^\phi N$ is $\pi$-separable. Hence, $H^\phi =H^{xy}$ for some $y\in H^\phi N$ by Lemma~\ref{Separ}.
 \qed\medskip

\section{Proof of Theorem~\ref{Th1} and Corollary~\ref{Cor}}
We divide our proof of Theorem~\ref{Th1} onto three parts. In the first part (Proposition~\ref{Determination} in Section \ref{1part}), we prove that if $H$ is a $\pi$-submaximal subgroup of a minimal simple group~$S$, then  $H$ can be found in that of Tables~\ref{L_2(2),2_notin_pi}--\ref{L_3(3),2,3} which corresponds to given $S$ and $\pi$. In the second part (Proposition~\ref{ProofOfSubmaximality} in Section \ref{2part}), we prove that every  $H$ in Tables~\ref{L_2(2),2_notin_pi}--\ref{L_3(3),2,3} are $\pi$-submaximal for corresponding $S$ and prove that the information on the $\pi$-maximality, conjugacy, intravariancy and the action of $\Aut(S)$ on the  conjugacy classes for this subgroups is true. Finally, in the third  part (Proposition~\ref{Pronormality} in Section \ref{3part}) we prove the pronormality of the $\pi$-submaximal subgroups in the minimal simple groups and, as a consequence, in the minimal nonsolvable groups (Corollary~\ref{Cor}).

\subsection{
The classification of $\pi$-submaximal subgroups\\ in minimal simple groups}\label{1part}

In this section, $S$ is a group in the Thompson list $\cal T$, that is, $S$ is  one of the following groups:
\begin{itemize}
\item[$(1)$] $L_2(q)$ where $q=2^p$, $p$~is a prime;

\item[$(2)$] $L_2(q)$ where where $q=3^p$,  $p$~is an odd prime;

 \item[$(3)$] $L_2(q)$ where $q$~is a prime such that $q>3$ and $q^2+1\equiv 0\pmod 5$;

  \item[$(4)$] $Sz(q)$ where  $q=2^p$, $p$~is an odd prime;

  \item[$(5)$] $L_3(3)$.
 \end{itemize}
  \noindent We identify $S$ with $\Inn(S)\cong S$. Let $\pi$ be a set of primes such that $|{\pi\cap\pi(S)}|>1$ and $\pi(S)\nsubseteq\pi$. The following statement gives a classification of $\pi$-submaximal subgroups in the minimal simple groups.

\begin{Pro}\label{Determination} {\it If $H\in\sm_\pi(S)$ where $S\in\cal T$, then $H$ appears in the corresponding column in that of  Tables {\rm \ref{L_2(2),2_notin_pi}--\ref{L_3(3),2,3}} which corresponds to~$S$ and $\pi$.}
\end{Pro}
\noindent{\sc Proof}. Let $H\in\sm_\pi(S)$. Lemma~\ref{SubmaxXSimple} implies that $H=K\cap S$ for some $K\in  \m_\pi(G)$ where $S\leq G\leq \Aut(S)$ and $G=KS$. By Lemma~\ref{WieHartSimple}, $N_S(H)/H$ is a $\pi'$-group and, in particular, $H\ne1$. Since $S$ is not a $\pi$-group, we have $K<G$ and $K\leq M$ for some maximal subgroup $M$ of~$G$. Note that $M\cap S<S$ (in fact, if $M\cap S=S$ it would be $G=KS\leq MS\leq M$, which contradicts $M\lessdot G$), so $M$ is solvable. Moreover, $K\in\m_\pi(M)=\Hall_\pi(M)$ in view of solvability of $M$. Hence by Lemma~\ref{HallSubgroup}, $$H=K\cap S=K\cap (M\cap S)\in\Hall_\pi(M\cap S).$$
It follows from Lemma~\ref{MaxSubgrAut} that one of the following cases holds:
\begin{itemize}
  \item[(I)]  $M\cap S\lessdot S$;
  \item[(II)]  $S\cong L_2(7)$, $G=\Aut(S)$, $M\cong D_{12}$ and $M\cap S\cong D_6$;
  \item[(III)]  $S\cong L_2(7)$, $G=\Aut(S)$, $M\cong D_{16}$ and $M\cap S\cong D_8$;
  \item[(IV)]  $S\cong L_3(3)$, $G=\Aut(S)$, $M\cong GL_2(3):2$ and $M\cap S\cong GL_2(3)$;
  \item[(V)]  $S\cong L_3(3)$, $G=\Aut(S)$, $M\cong 3_+^{1+2}:D_8$ and $M\cap S\cong 3_+^{1+2}:C_2^2$.
\end{itemize}

 Firstly, we consider Cases (II)--(V) keeping in mind that $1<H\in\Hall_\pi(M\cap S)$.
 \medskip

 Assume that Case (II) holds. Then $S=L_2(q)$, where $q=7$ and $M\cap S\cong D_6$.

 If $2\notin\pi$, then $|H|=3$ and $H$  appears  in the 2-nd row of Table~\ref{L_2(q),2_notin_pi} for $\pi_+=\{3\}$.

 If $3\notin\pi$, then $|H|=2$ and $2\in \pi$.  In this case, $H$ is contained in a Sylow 2-subgroup $P$ of $S$. But $1<N_P(H)/H \leq N_S(H)/H$, so $N_S(H)/H$ is not a $\pi'$-group, which  contradicts  Lemma~\ref{WieHartSimple}. Thus, this case is impossible.

 If $2,3\in \pi$, then $H=M\cap S$ and $H$  appears  in the 2-nd row of  Table~\ref{L_2(q),2_in_pi} for $\pi_+\ne\{2\}$.

 \medskip

Assume that  Case (III) holds. Then $S=L_2(q)$, where $q=7$ and $M\cap S\cong D_8$. In this case, $H\in\Hall_\pi(M\cap S)$ implies that $2\in\pi$,  $H=M\cap S$ is a Sylow 2-subgroup of $S$, and $H$  appears  in the 3-rd row of  Table~\ref{L_2(q),2_in_pi} for $q\equiv -1\pmod4$.
\medskip

Assume that  Case (IV) holds. Then $S=L_3(3)$ and $M\cap S\cong GL_2(3)$.

 If $2\notin\pi$, then $|H|=3$ and $H$ is contained in a Sylow 3-subgroup $P$ of $S$. But $1<N_P(H)/H\leq N_S(H)/H$, so $N_S(H)/H$ is not a $\pi'$-group, which contradicts Lemma~\ref{WieHartSimple}. Thus, this case is impossible.

 If $3\notin\pi$, then  $H\cong SD_{16}$, $H$ coincides with a Sylow 2-subgroup of $S$, and $\pi\cap\pi(S)=\{2,13\}$. Hence,   $H$  appears  in the 2-nd row of  Table~\ref{L_3(3),2,13}.

 If $2,3\in \pi$, then $H=M\cap S$ and $H$  appears  in the 3-rd row of  Table~\ref{L_3(3),2,3} for $\pi\cap\pi(S)=\{2,3\}$.

 \medskip

Assume that  Case (V) holds. Then $S=L_3(3)$ and $M\cap S\cong 3_+^{1+2}:C_2^2$.

 If $2\notin\pi$, then $\pi\cap\pi(S)=\{3,13\}$, $H\cong 3_+^{1+2}$ and $H$ coincides with a Sylow 3-subgroup of $S$. In this case,  $H$  appears  in the 2-nd row of  Table~\ref{L_3(3),3,13}.

 If $3\notin\pi$, then  $|H|=4$  and $2\in \pi$.  In this case, $H$ is contained in a Sylow 2-subgroup $P$ of $S$. But $1<N_P(H)/H\leq N_S(H)/H$. Hence $N_S(H)/H$ is not a $\pi'$-group, contrary to Lemma~\ref{WieHartSimple}. Thus, this case is impossible.

 If $2,3\in \pi$, then $H=M\cap S$ and $H$  appears  in the 2-nd row of  Table~\ref{L_3(3),2,3} for $\pi\cap\pi(S)=\{2,3\}$.

 \medskip

 Now we consider Case (I). In this case, $H$ coincides with some nontrivial $\pi$-Hall subgroup of a maximal subgroup $U=M\cap S$ of $S\in \cal T$. By using Lemma~\ref{MaxSubgrSimple}, we consider the nontrivial $\pi$-Hall subgroups of maximal subgroups $U$ of all groups  $S\in \cal T$.
  \medskip

 Case (I)(1): $S=L_2(q)$ where $q=2^p$, $p$~is a prime. In this case, the numbers $2$, ${q-1}$ and ${q+1}$ are pairwise coprime and hence $$\pi(S)=\{2\}\,\dot{\cup}\,\pi(q-1)\,\dot{\cup}\,\pi(q-1)$$ (the dot over the symbol of union means that we have a union of pairwise disjoint sets). Let
 $$\pi_\varepsilon=\pi{\cap}\pi(q-\varepsilon)\text{ for }\varepsilon\in\{+,-\}.$$
 By Lemma~\ref{MaxSubgrSimple}, $U$ is one of the following groups: either
     $E_q:C_{(q-1)}$, or
  $D_{2(q-\varepsilon)}$, where $\varepsilon\in\{+,-\}$.

  If $2\notin\pi$, then $$\pi\cap\pi(S)=\pi_+\,\dot{\cup}\,\pi_-$$ and $H$ must be
 a nontrivial cyclic group of order $(q-\varepsilon)_\pi$ and so $\pi_\varepsilon\ne\varnothing$. Thus, $H$  appears  in Table~\ref{L_2(2),2_notin_pi}.

 If $2\in\pi$, then $$\pi\cap\pi(S)=\{2\}\,\dot{\cup}\,\pi_+\,\dot{\cup}\,\pi_-$$ and $H$ must be either the Frobenius group in the 1-st row of Table~\ref{L_2(2),2_in_pi} or the dihedral group  $D_{2(q-\varepsilon)_\pi}$, where $\varepsilon\in\{+,-\}$. Moreover, in the last case,   $\pi_\varepsilon\ne\varnothing$  since otherwise $|H|=2$ and $H$ is contained as a proper subgroup in a Sylow 2-subgroup $P$ of $S$, contrary to Lemma~\ref{WieHartSimple}. Thus $H$  appears  in Table~\ref{L_2(2),2_in_pi}.

 \medskip

Case (I)(2): $S=L_2(q)$ where $q=3^p$, $p$~is an odd prime. In this case, $$\pi(S)=\{3\}\cup\pi(q-1)\cup\pi(q-1).$$ Moreover,  $3$ is coprime with both ${q-1}$ and ${q+1}$ and $(q-1, q+1)=2$. As in above, we set
 $$\pi_\varepsilon=\pi\cap\pi(q-\varepsilon)\text{ for }\varepsilon\in\{+,-\}.$$

 By Lemma~\ref{MaxSubgrSimple}, $U$ is one of the following groups: either
     $E_q:C_{\frac{1}{2}(q-1)}$, or
  $D_{q-\varepsilon}$, where $\varepsilon\in\{+,-\}$, or $A_4$.

  If $2,3\notin\pi$, then $$\pi\cap\pi(S)=\pi_+\,\dot{\cup}\,\pi_-$$ and $H$ must be
 a nontrivial cyclic group of order $(q-\varepsilon)_\pi$, and so $\pi_\varepsilon\ne\varnothing$. Thus, $H$  appears  in Table~\ref{L_2(3),2_notin_pi}.

 If $2\notin\pi$ and $3\in\pi$, then $$\pi\cap\pi(S)=\{3\}\,\dot{\cup}\,\pi_+\,\dot{\cup}\,\pi_-$$ and $H$ must be  either the Frobenius group in the 1-st row of Table~\ref{L_2(3),2_notin_pi} or
 a nontrivial cyclic group of order $(q-\varepsilon)_\pi$ and hence $\pi_\varepsilon\ne\varnothing$.

Assume that $\varepsilon=+$ and $H$ is  cyclic  of order $(q-1)_\pi$.  We claim that $H\notin\sm_\pi(G)$. Indeed, if $H\in\sm_\pi(G)$, then $H=K\cap S$ for some $K\in\m_\pi(\Aut(S))$ by Lemma~\ref{SubmaxXSimple}. It is well-known that ${|\Aut(S)/S|=2p}$ (see~\cite[table~8.1]{Bray}, for example). Consider two cases: $p\notin\pi$ and $p\in\pi$.

 If $p\notin\pi$, then $K\leq S$ and $H=K$. In particular, in this case, $H\in\m_\pi(S)$. But it is easy to see that $H\in\Hall_{\pi_+}(S)$ and  $S\in\D_{\pi_+}$ by Lemma~\ref{NilpHall}. A maximal subgroup of the type $E_q:C_{\frac{1}{2}(q-1)}$ of $S$ contains some $\pi_+$-Hall subgroup of $S$ and $H$ is conjugate to this subgroup. It means that $H$ normalizes a Sylow 3-subgroup of $S$, contrary to $H\in\m_\pi(S)$.

 If $p\in\pi,$ then, in view of $(p,q+1)=(p,3^p+1)=1$, we have by Lemma~\ref{DpiInL2q} that
 $$\Aut(S)\in\D_\tau\text{ where } \tau=\pi_+\cup\{3,p\}.$$
 Since $K\in \m_\pi(\Aut(S))$ and $K$ is a $\tau$-group, we obtain that $$K\in\m_\tau(\Aut(S)) =\Hall_\tau(\Aut(S)).$$ But $|K|=p|H|$ is not divisible by 3, a contradiction again.

 Thus, if  $2\notin\pi$ and $3\in\pi,$ then $H$   is not a  cyclic group  of order $(q-1)_\pi$. In the remaining cases, $H$  appears  in Table~\ref{L_2(3),2_notin_pi}.

 Suppose that $2\in\pi$. Then the maximal subgroups containing a Sylow 2-subgroup of $S$ are $D_{q+1}$ and $A_4$. Hence $H$ is a 2-group if and only if $\pi_-=\{2\}$ and $H=D_{(q+1)_2}=D_{(q+1)_\pi}$. In particular, such $H$ can not be contained  in $E_q:C_{\frac{1}{2}(q-1)}$ and
  $D_{q-1}$. Now it is easy to see that if $3\notin\pi$, then $H$ coincides with a $\pi$-Hall subgroup of $D_{q-\varepsilon}$, where $\varepsilon\in\{+,-\}$; and if $3\in\pi$, then  $H$ can not coincide with a Sylow 2-subgroup $P$ by Lemma~\ref{WieHartSimple} in view of $N_S(P)\cong A_4$. Hence $H$ must coincides with one of the following groups $E_q:C_{\frac{1}{2}(q-1)_\pi}$, $D_{(q-1)_\pi}$, $D_{(q+1)_\pi}$ or~$A_4$. In all these cases, $H$  appears  in Table~\ref{L_2(3),2_in_pi}.

\medskip

Case (I)(3): $S=L_2(q),$ where $q$~is a prime such that $q>3$ and $q^2+1\equiv 0\pmod 5$. In this case, $$\pi(S)=\{q\}\cup\pi(q-1)\cup\pi(q-1),$$   $q$ is coprime with both ${q-1}$ and ${q+1}$ and $(q-1, q+1)=2$. Let
 $$\pi_\varepsilon=\pi\cap\pi(q-\varepsilon)\text{ for }\varepsilon\in\{+,-\}$$ as usual.

 In this case, $\Aut(S)\cong PGL_2(q)$. By Lemma~\ref{MaxSubgrSimple}, $U$ is one of the following groups: or
     $C_q:C_{\frac{1}{2}(q-1)}$, or
  $D_{q-\varepsilon}$ where $\varepsilon\in\{+,-\}$, or one of $A_4$ and $S_4$  for $q\equiv\pm3\pmod 8$ and  $q\equiv\pm1\pmod 8$, respectively. It follows from \cite[table~8.1]{Bray} that, with the exception of the case~$U\cong S_4$, we have $U=V\cap S$ for some $V\lessdot\Aut(S)$ such that $\Aut(S)=SV$ and $|V:U|=2$.

  Suppose that $2\notin\pi$. Since $|\Aut(S):S|=2$, $H$ must be a $\pi$-maximal subgroup of $S$ by Lemma~\ref{SubmaxXSimple}. If $H$ is contained in one of $A_4$ or $S_4$, then $|H|=3$ and, since any Sylow 3-subgroup of $S$ is cyclic, $H$ is contained in $D_{q-\varepsilon}$ with $q\equiv\varepsilon\pmod 3$. Moreover, if $q\in \pi$, then the unique $\pi$-Hall subgroup of $D_{q-1}$ is not $\pi$-maximal in $S$. Indeed, this subgroup is a cyclic $\pi_+$-Hall subgroup of $S$.  Lemma~\ref{NilpHall} implies that this subgroup is conjugate to a subgroup in some Frobenius group $C_q:C_{\frac{1}{2}(q-1)}\leq S$ and normalizes a Sylow $q$-subgroup of~$S$. Hence, if $2\notin\pi$, then $H$ appears  in Table~\ref{L_2(q),2_notin_pi}.

  Suppose that $2\in\pi$. In order to show that $H$ appears  in Table~\ref{L_2(q),2_in_pi}, it is sufficient to show that a $\pi$-Hall subgroup $H$ of $D_{q-\varepsilon}$ is not $\pi$-submaximal in $S$ if $\pi_\varepsilon=\{2\}$, $3\in\pi$ and $q\not\equiv \varepsilon\pmod8$. Indeed, if $q\equiv-\varepsilon\pmod8$, then $|H|=(q-\varepsilon)_2=2<|S|_2$ and $H<P$ for some Sylow 2-subgroup  $P$ of $S$, which contradicts Lemma~\ref{WieHartSimple}; if $q\equiv\pm3\pmod8$, then $H$ is an elementary  abelian Sylow 2-subgroup of $S$, and $|H|=4$, which contradicts Lemma~\ref{WieHartSimple} since, by the Sylow theorem, $H$ is contained as a normal subgroup in a maximal subgroup of $S$ isomorphic to $A_4$ and this subgroups is a $\pi$-group.

  \medskip

  Case (I)(4): $S=Sz(q)$ where $q=2^p$, $p$~is an odd prime. In this case, $$\pi(S)=\{2\}\,\dot{\cup}\,\pi(q-1)\,\dot{\cup}\,\pi(q-r+1)\,\dot{\cup}\,\pi(q+r+1)$$ where $r=\sqrt{2q}=2^{(p+1)/2}$.   Let
 $$\pi_0=\pi\cap \pi(q-1)\text{ and }\pi_\varepsilon=\pi\cap\pi(q-\varepsilon r+1)\text{ for }\varepsilon\in\{+,-\}.$$

   By Lemma~\ref{MaxSubgrSimple}, $U$ is one of the following groups: or  ${E_q^{1+1}:C_{(q-1)}}$,  or
 $C_{(q-1)_\pi}:C_4$, or
 $D_{2(q- \varepsilon r+1)_\pi}$ where $\varepsilon\in\{+,-\}$.

 Note that $S$ has exactly one conjugacy class of cyclic subgroups of order $q-1$ (see the character table of $S$ in \cite[theorem~13]{SuzCl}). Hence, if $2\notin\pi$, then any $\pi$-Hall subgroup of $E_q^{1+1}:C_{(q-1)}$ is contained as a $\pi$-Hall subgroup in some maximal subgroup of kind  ${C_{(q-1)_\pi}:C_4}$. Thus, in the case when $2\notin\pi$, $H$  appears  in Table~\ref{Sz,2_notin_pi}.

Now assume that $2\in\pi$. If one of the sets $\pi_0$, $\pi_+$, $\pi_-$ is empty, then any $\pi$-Hall subgroup of the respective subgroups: or ${C_{(q-1)_\pi}:C_4}$, or
 $D_{2(q- r+1)}$, or $D_{2(q- r+1)}$,  is a $2$-group, but is not a Sylow 2-subgroup. Then Lemma~\ref{WieHartSimple} implies that case where $H$ is a $\pi$-Hall subgroup in one of these maximal subgroup of $S$ is impossible. In the remaining cases,   $H$  appears  in Table~\ref{Sz,2_in_pi}.

  \medskip

  Case (I)(5): $S=L_3(3)$. In this case, $\pi(S)=\{2,3,13\}.$ By Lemma~\ref{MaxSubgrSimple}, $U$ is one of the following groups:
or $E_{3^2}:GL_2(3)$, or
  $C_{13}:C_{3}$, or
 $S_4$.   Since $\pi(S)\nsubseteq\pi$ and $|\pi\cap\pi(S)|>1$, the intersection $\pi\cap\pi(S)$ coincides with one of sets  $\{3,13\}$,   $\{2,13\}$, and $\{2,3\}$.

 Suppose that $\pi\cap\pi(S)=\{3,13\}$. Then a $\pi$-Hall subgroup $H$ of  $E_{3^2}:GL_2(3)$ coincides with a Sylow 3-subgroup $3_+^{1+2}$ of $S$ and appears in  Table~\ref{L_3(3),3,13}. A $\pi$-Hall subgroup $H$ of $C_{13}:C_{3}$ is $C_{13}:C_{3}$ itself  and $H$ appears in  Table~\ref{L_3(3),3,13}.  A $\pi$-Hall subgroup $H$ of $S_4$ is of order 3 and Lemma~\ref{WieHartSimple} implies that this case is impossible.

 Suppose that $\pi\cap\pi(S)=\{2,13\}$. Then a $\pi$-Hall subgroup $H$ of  $E_{3^2}:GL_2(3)$ coincides with a semi-dihedral Sylow 2-subgroup $SD_{16}$ of both $E_{3^2}:GL_2(3)$ and $S$, and $H$ appears in  Table~\ref{L_3(3),2,13}. A $\pi$-Hall subgroup $H$ of $C_{13}:C_{3}$ is cyclic of order 13 and  appears in  Table~\ref{L_3(3),2,13}.   A $\pi$-Hall subgroup $H$ of $S_4$ is of order 8 and is not a Sylow 2-subgroup of $S$. We may exclude this case by Lemma~\ref{WieHartSimple}.

 Finally, suppose that $\pi\cap\pi(S)=\{2,3\}$. Then a $\pi$-Hall subgroup $H$ of  $E_{3^2}:GL_2(3)$ is  $E_{3^2}:GL_2(3)$ itself and $H$ appears in  Table~\ref{L_3(3),2,3}. A $\pi$-Hall subgroup $H$ of $C_{13}:C_{3}$ is cyclic of order 3 and  Lemma~\ref{WieHartSimple} implies that this case is impossible.  A $\pi$-Hall subgroup $H$ of $S_4$ is $S_4$ itself and  it appears in  Table~\ref{L_3(3),2,3}.
\qed\medskip

\subsection{ 
The $\pi$-submaximality, $\pi$-maximality, conjugacy,\\ and intravariancy for subgroups in Tables~\ref{L_2(2),2_notin_pi}--\ref{L_3(3),2,13}}\label{2part}

In this section, as in above, $S$ is a group in the Thompson list $\cal T$.

We prove that the converse of Proposition~\ref{Determination} holds, that is, the subgroups placed in Tables~\ref{L_2(2),2_notin_pi}--\ref{L_3(3),2,13} are $\pi$-submaximal under corresponding conditions. Moreover, we prove that the information about the $\pi$-maximality, the conjugacy, the action of $\Aut(S)$ on the conjugacy classes of this subgroups, and the intravariancy in this tables is correct.

\begin{Pro}\label{ProofOfSubmaximality} {\it Let $S\in \cal T$,  $\pi$ be a set of primes such that $|{\pi\cap\pi(S)}|>1$ and ${\pi(S)\nsubseteq\pi}$, and let $H$ be a subgroup of $S$ placed in corresponding column of one of Tables {\rm \ref{L_2(2),2_notin_pi}--\ref{L_3(3),2,3}} which corresponds to~$S$. Then
\begin{itemize}
  \item[$(1)$] $H\in\sm_\pi(S)$.
  \item[$(2)$] $H\in\Hall_\pi(M)$ for every $M\lessdot S$ such that $H\leq M$ (and so $H\in\m_\pi(S)$) with the exception of the following cases:
  \begin{itemize}
    \item[{\rm (2a)}] $S=L_2(q)$ for some prime $q$ such that $q>3$ and $q^2+1\equiv 0\pmod 5$;\\ $2,3\in \pi$, $H=D_{(q-\varepsilon)_\pi}$, $\varepsilon\in\{+,-\}$,
       $\pi\cap\pi(q-\varepsilon)=\{2\}$ and $$q\equiv-\varepsilon7
       \pmod{48}.$$
    In this case, $H$ is a Sylow $2$-subgroup of $S$, $H\cong D_8$,  and  $H<S_4\lessdot S$.
    \item[{\rm (2b)}] $S=L_2(q)$ for some prime $q$ such that $q>3$ and $q^2+1\equiv 0\pmod 5$;\\ ${2,3\in \pi}$, ${H=D_{(q-\varepsilon)_\pi}}$, ${\varepsilon\in\{+,-\}}$,
              ${{\pi\cap\pi(q-\varepsilon)}=\{2,3\}}$ and $${q\equiv\varepsilon7,\varepsilon31\pmod{72}}.$$
    In this case, $H\cong S_3\cong D_6$, $H$ is a $\{2,3\}$-Hall subgroup of the normalizer in $S$ of a Sylow $3$-subgroup and $H<S_4\lessdot S$.
    \item[{\rm (2c)}] $S=L_3(3)$, $\pi\cap\pi(S)=\{2,3\}$, and ${H={3_+^{1+2}:C_2^2}}$. In this case, ${H<{E_{3^2}:GL_2(3)}\lessdot S}$.
    \item[{\rm (2d)}] $S=L_3(3)$, $\pi\cap\pi(S)=\{2,3\}$, and $H=GL_2(3)$. In this case, ${H<{E_{3^2}:GL_2(3)}\lessdot S}$.
  \end{itemize}
  \item[$(3)$] $S$ has a unique conjugacy class of $\pi$-submaximal subgroups isomorphic to $H$ with the exception of the following cases when $S$ has exactly two conjugacy classes of $\pi$-submaximal subgroups isomorphic to $H$ that are fused in $\Aut(S)$:
  \begin{itemize}
    \item[{\rm (3a)}] $S=L_2(q)$ where  $q$~is a prime such that $q>3$, $q^2+1\equiv 0\pmod 5$, $q\equiv\pm1\pmod 8$, $2,3\in \pi$, and $H=S_4$;
    \item[{\rm (3b)}] $S=L_3(3)$, $\pi\cap\pi(S)=\{2,3\}$, and $H=E_{3^2}:GL_2(3)$.
  \end{itemize}
  \item[$(4)$] $H$ is intravariant in $S$ excepting the cases determined in $(3)$ where $S$ has two conjugacy classes of $\pi$-submaximal subgroups isomorphic to~$H$.
\end{itemize}}
\end{Pro}
\noindent{\sc Proof}.
Firstly, we prove~(2).

Non-$\pi$-maximality of $H$ and corresponding inclusions in cases (2b) or (2c) follow from~\cite{Atlas}.

Suppose case (2a) holds. Then  $S=L_2(q)$, $2,3\in \pi$, $H=D_{(q-\varepsilon)_\pi}$ where  $q$~is a prime such that $q>3$ and $q^2+1\equiv 0\pmod 5$, $\varepsilon\in\{+,-\}$.

Consider the case where $\pi\cap\pi(q-\varepsilon)=\{2\}$ and $q\equiv-\varepsilon7
\pmod{48}$. Then $q\equiv \varepsilon\pmod8$. This means that $H$ coincides with a Sylow 2-subgroup of $S$ and, moreover, $S$ contains a maximal subgroup isomorphic to~$S_4$ by Lemma~\ref{MaxSubgrSimple}. Since $q\equiv -7\varepsilon\pmod{16}$, we have $|H|=(q-\varepsilon)_\pi=8=|S_4|_2$. The Sylow theorem implies that $H$ is conjugate to a Sylow 2-subgroup of $S_4$. Thus, $H$ is not $\pi$-maximal in $S$ in view of $2,3\in\pi$ and $S_4$ is a $\pi$-group.

Now consider the case where $\pi\cap\pi(q-\varepsilon)=\{2,3\}$ and $q\equiv\varepsilon7,\varepsilon31\pmod{72}$. Then $q\equiv -\varepsilon\pmod8$. It means that $S$ contains a maximal subgroup $M$ isomorphic to~$S_4$ by Lemma~\ref{MaxSubgrSimple}. But it is easy to calculate that $|H|=(q-\varepsilon)_\pi=6$ and so $H\cong S_3$. Moreover, $H$ contains a Sylow 3-subgroup $P$ of $S$ since $|S|_3=3$. It follows that $H$ is contained in $N_S(P)$. By considering the maximal subgroups of $S$ given in Lemma~\ref{MaxSubgrSimple}, it is easy to see that $N_S(P)=D_{q-\varepsilon}$, so $H$ is a $\{2,3\}$-Hall subgroup of  $N_S(P)$.  By the Sylow theorem, we can assume that $P<M$. Since $N_M(P)\cong S_3$, we obtain $N_M(P)\in\Hall_{\{2,3\}}(N_S(P))$. The Hall theorem and the solvability of $N_S(P)$ imply that $H$ is conjugate with  $N_M(P)$ and $H$ is not $\pi$-maximal in $S$ since $2,3\in\pi$ and $S_4$ is a $\pi$-group.

In order to complete the proof of (2), it is sufficient to show that, in the remaining cases, $H$ is a $\pi$-Hall subgroup in every maximal subgroup $M$ of $S$ containing $H$. By  Lemma~\ref{MaxSubgrSimple}, an easy calculation show that one of the following cases holds:
\begin{itemize}
    \item[{\rm (I)}] up to isomorphism, there is a unique maximal subgroup $M$ of $S$ such that $|H|$ divides $|M|$, $H\leq M$ and $H$ is a $\pi$-Hall subgroup of $M$.
    \item[{\rm (II)}] $S=L_2(q)$ where $q=2^p$ for a prime $p$, $2\notin\pi$, $H=C_{(q-1)_\pi}$, and, up to isomorphism, there are exactly two maximal subgroup $M$ of $S$ such that $|H|$ divides $|M|$, namely either $M=E_q:C_{(q-1)}$ or
  $M=D_{2(q-1)}$. In this case, $H$ is a $\pi$-Hall subgroup of every maximal subgroup  containing $H$.
  \item[{\rm (III)}] $S=L_2(q)$ where $q=2^p$ for a prime $p$,  $2\in\pi$, $H=D_{2(q-1)_\pi}$, and, up to isomorphism, there are exactly two maximal subgroup $M$ of $S$ such that $|H|$ divides $|M|$, namely either $M=E_q:C_{(q-1)}$ or
  $M=D_{2(q-1)}$. In this case $H$ is not isomorphic to a subgroup of $E_q:C_{(q-1)}$ and is contained only in $M=D_{2(q-1)}$. Hence $H$ is a $\pi$-Hall subgroup of every maximal subgroup  containing $H$.
    \item[{\rm (IV)}] $S=L_2(q)$ where either $q=3^p$ for odd prime or $q>3$ is a prime, $2\notin\pi$,  $H=C_{(q-1)_\pi}$ and, up to isomorphism, there are exactly two maximal subgroup $M$ of $S$ such that $|H|$ divides $|M|$, namely either $M=E_q:C_{\frac{1}{2}(q-1)}$ or
  $M=D_{q-1}$.  In this case, $H$ is a $\pi$-Hall subgroup of every maximal subgroup  containing $H$.
    \item[{\rm (V)}] $S=Sz(q)$ where $q=2^p$ for a prime $p$, $2\notin\pi$, $H=C_{(q-1)_\pi}$, and, up to isomorphism, there are exactly two maximal subgroup $M$ of $S$ such that $|H|$ divides $|M|$, namely either $M=E^{1+1}_q:C_{(q-1)}$ or
  $M=C_{q-1}:C_4$.  In this case, $H$ is a $\pi$-Hall subgroup of every maximal subgroup  containing $H$.
   \item[{\rm (VI)}] $S=Sz(q)$ where $q=2^p$ for a prime $p$,  $2\in\pi$, $H=C_{(q-1)_\pi}:C_4$, and, up to isomorphism, there are exactly two maximal subgroup $M$ of $S$ such that $|H|$ divides $|M|$, namely either $M=E^{1+1}_q:C_{(q-1)}$ or
  $M=C_{q-1}:C_4$. In this case $H$ is not isomorphic to a subgroup of $E^{1+1}_q:C_{(q-1)}$ and is contained only in $M=C_{q-1}:C_4$. Hence $H$ is a $\pi$-Hall subgroup of every maximal subgroup  containing $H$.
  \item[{\rm (VII)}]  $S=L_2(q)$ where  $q$~is a prime such that $q>3$ and $q^2+1\equiv 0\pmod 5$, $2\in \pi$, $H=D_{(q-\varepsilon)_\pi}$ for some $\varepsilon\in\{+,-\}$,  $|H|$ divides $24$, and, up to isomorphism, there are exactly two maximal subgroup $M$ of $S$ such that $|H|$ divides $|M|$, namely either $M=D_{q-\varepsilon}$ or
  $M\in\{A_4,S_4\}$. This case will be separately  considered  below.
  \item[{\rm (VIII)}]  $S=L_2(q)$ where  $q$~is a prime such that $q>3$ and $q^2+1\equiv 0\pmod 5$, $2,3\in \pi$, $H\in\{A_4,S_4\}$, $|H|$ divides $q-\varepsilon$  for some $\varepsilon\in\{+,-\}$ and, up to isomorphism, there are exactly two maximal subgroup $M$ of $S$ such that $|H|$ divides $|M|$, namely either $M=D_{q-\varepsilon}$ or
  $M\in\{A_4,S_4\}$. In this case $H$ is not isomorphic to any subgroup of  $D_{q-\varepsilon}$, so $H$ itself is maximal and coincide with its $\pi$-Hall subgroup.
  \item[{\rm (IX)}]  ${S=L_3(3)}$, ${\pi=\{2,3\}}$, ${H=S_4}$, and, up to isomorphism, there are exactly two maximal subgroup $M$ of $S$ such that $|H|$ divides $|M|$, namely either ${M=E_{3^2}:GL_2(3)}$ or
  $M=S_4$. In this case, $H$ itself is maximal and  coincide with its $\pi$-Hall subgroup.

      \end{itemize}
   Thus, in order to complete the proof of (2), we need to consider  Case~(VII):  $S=L_2(q)$ where  $q$~is a prime such that $q>3$ and $q^2+1\equiv 0\pmod 5$, $2\in \pi$, $H=D_{(q-\varepsilon)_\pi}$ for some $\varepsilon\in\{+,-\}$,  $|H|$ divides $24$, and, up to isomorphism, there are exactly two maximal subgroup $M$ of $S$ such that $|H|$ divides $|M|$, namely either $M=D_{q-\varepsilon}$ or
  $M\in\{A_4,S_4\}$. It is clear that $H$ is a $\pi$-Hall subgroup of $D_{q-\varepsilon}$.  We need to show that either $H$ is not contained in $M\in\{A_4,S_4\}$, or $H$ is a $\pi$-Hall subgroup in such~$M$, or one of the exceptional cases (2a) and (2b) holds.

  Note that we have the condition in Table~\ref{L_2(q),2_in_pi} for $H$ that ether $\pi_\varepsilon\ne\{2\}$, or $3\notin\pi$, or $q\equiv \varepsilon\pmod 8$.

   Suppose that $\pi_\varepsilon\ne\{2\}.$ Since $|H|=(q-\varepsilon)_\pi$ divides 24, the assumption  $\pi_\varepsilon\ne\{2\}$ means that $3\in\pi$, $3$ divides $H$, and $M\in\{A_4,S_4\}$ is a $\pi$-group. Since $H$ must contain a dihedral group of order 6 in this case, $M\ne A_4$ and so $M=S_4$. This implies that $q\equiv\pm1\pmod8$. But $S_4$ does not contained a dihedral group of order $12$ and we have that $|H|=6$. In particular, $$2=|H|_2= (q-\varepsilon)_2$$ and  $q\equiv -\varepsilon\pmod 4$. Thus, $q\equiv -\varepsilon\pmod 8$. Moreover, $|H|$ divides 24 implies that $q-\varepsilon$ is not divisible by $9$ and hence $q\equiv \pm3+\varepsilon\pmod 9$. It is not difficult to prove that $${q\equiv -\varepsilon\!\!\pmod 8}\,\,\text{ and }\,\,{q\equiv \pm3+\varepsilon\!\!\pmod 9}\,\,\,
  \text{ if and only if }\,\,\, {q\equiv \varepsilon7, \varepsilon31\!\!\pmod{72}}.$$ Thus,  exceptional case (2b) holds.

   Thus, we can consider that $\pi_\varepsilon=\{2\}$. Suppose that $3\notin\pi$. Since $|H|$ divides $24$, we have $H=D_8$ or $H=D_4=C_2\times C_2$. In the first case, $M=S_4$ and $H$ is a $\pi$-Hall subgroup of both $D_{q-\varepsilon}$ and $S_4$, and in the second case, $M=A_4$ and $H$ is $\pi$-Hall in both $D_{q-\varepsilon}$ and $A_4$.

   Consider the last case when $\pi_\varepsilon=\{2\}$ and $3\in\pi$. Then $q\equiv \varepsilon\pmod 8$ and condition $|H|$ divides $24$ implies that $|H|=8$. Hence $q\not\equiv \varepsilon\pmod {16}$. This implies that  $q\equiv 8+\varepsilon\pmod {16}$. Moreover, $q\equiv -\varepsilon\pmod 3$ since otherwise $3\in \pi_\varepsilon=\{2\}$. It is not difficult to prove that $$q\equiv 8+\varepsilon\!\!\pmod {16}\,\,\,\text{ and }\,\,\,q\equiv -\varepsilon\!\!\pmod 3\,\,\,\text{ if and only if }\,\,\,q\equiv-\varepsilon7\!\!
   \pmod{48}.$$ Thus, the  exceptional case (2a) holds.

This completes the proof of~(2).

\medskip

Now we prove (1). In view of (2), it is sufficient to prove that if one of the exceptional cases (2a)--(2d) in (2) holds, then $H=K\cap S$ for some $K\in\m_\pi(G)$ where $G=\Aut(S)$.

Suppose one of cases (2a)--(2b) holds. Then $S=L_2(q)$, $2,3\in \pi$, $H=D_{(q-\varepsilon)_\pi}$ where  $q$~is a prime such that $q>3$, $q^2+1\equiv 0\pmod 5$, $\varepsilon\in\{+,-\}$, and $G=PGL_2(q)$. It follows from~\cite[table~8.1]{Bray} that $H\leq U\lessdot G$ where $U\cong D_{2(q-\varepsilon)}$. Let $K\cong D_{2(q-\varepsilon)_\pi}$ be a $\pi$-Hall subgroup of $U$. Then $G=SK$. Assume that $K\notin\m_\pi(G)$. Then $K<L$ for some $\pi$-maximal subgroup $L$ of $G$. It follows from  $\pi(S)\nsubseteq\pi$ that $S\nleq L$. Let $V$ be a maximal subgroup of $G$ such that    $L\leq V$.   In view of Lemma~\ref{MaxSubgrAut}, either $V\cap S\lessdot S$ or $q=7$ and $V\cap S\in\{D_6, D_8\}$. It is easy to see that the the both arithmetic conditions $q\equiv-\varepsilon7\pmod{48}$ in~(2a) and  $q\equiv\varepsilon7,\varepsilon31\pmod{72}$ in (2b) imply
that $q\equiv\pm1\pmod8$, so $S$ has no maximal subgroup isomorphic to $A_4$. Moreover, $VS=G=PGL_2(q)$. If we assume that $V\cap S\lessdot S$ and $V\cap S\cong S_4$, then the conjugacy class in $S$ of maximal subgroups isomorphic to $S_4$ would be invariant  under $G=\Aut(S)$. But it is not so, see Table~\ref{MaxL_2(q)}. Thus, $V\cap S$ is one of the following groups: $C_q:C_{\frac{1}{2}(q-1)}$, $D_{q+\varepsilon}$, and $D_{q-\varepsilon}$. Since $V\leq N_G(V\cap S)$, we have $V=N_G(V\cap S)$ and $V$ coincides with one of  $C_q:C_{q-1}$, $D_{2(q+\varepsilon)}$, and $D_{2(q-\varepsilon)}$. But $V$ contains  $K\cong D_{2(q-\varepsilon)_\pi}$ and $K\notin\Hall_\pi(V)$. We can exclude the case where $V=C_q:C_{q-1}$ since in this case every dihedral subgroup of $V$ is isomorphic to $D_{2q}$ but $((q-\varepsilon)_\pi,q)=1$. If $V= D_{2(q+\varepsilon)}$, then $(q-\varepsilon)_\pi$ divide $(q-\varepsilon,q+\varepsilon)=2$, but $(q-\varepsilon)_\pi=|H|\in\{6,8\}$; a contradiction. The last case when $V=D_{2(q-\
\varepsilon)}$ is impossible in view of $K\notin\Hall_\pi(V)$. Thus, $K\in\m_\pi(G)$ and $H=K\cap S\in\sm_\pi(S)$.

In the cases (2c) and (2d), $\pi\cap\pi(S)=\{2,3\}$. Let $G=\Aut(S)\cong L_3(3):C_2$.  Lemma~\ref{MaxSubgrAut} implies that $H=V\cap S$ where $V\lessdot G$ and $V=GL_2(3):C_2$ in the case (2c) and $V\cong 3_+^{1+2}:D_8$  in the case (2d). Thus $V\in\m_\pi(G)$ and $H\in\sm_\pi(S)$.
Therefore, (1) is proved.

\medskip

Prove (3). By Lemma~\ref{MaxSubgrSimple} (Tables~\ref{MaxL_2(q)} and~\ref{MaxL_3(3)}), if one of the exceptional cases (3a) or (3b) holds, then $H$ is a maximal subgroup of corresponding $S$, any maximal subgroup of $S$ which is not isomorphic to $H$ does not contain a subgroup isomorphic to~$H$, and $S$ has exactly two conjugacy classes of (maximal) subgroups  isomorphic to $H$ interchanged by every non-inner automorphism of $S$.

We need to prove that in the remaining cases, every $K\in\sm_\pi(S)$ isomorphic to $H$ is conjugate to $H$. By (2), if there is a maximal subgroup $M$ of $S$ containing both $H$ and $K$, then $H,K\in\Hall_\pi(M)$. The solvability of $M$ implies that $H$ and $K$ are conjugate in $M$ in view the Hall theorem. Hence we can consider that $H$ and $K$ are contained in non-conjugate maximal subgroups $M$ and $N$, respectively.

If $M$ and $N$ are isomorphic, then   Lemma~\ref{MaxSubgrSimple} implies that either (a) $S=L_2(q)$, $q>3$ is a prime, and $M\cong N\cong S_4$, or (b) $S=L_3(3)$ and $M\cong N\cong E_{3^2}:GL_2(3)$. Since $H\in\Hall_\pi(M)$ and $K\in\Hall_\pi(N)$ in view of (2) and $H\ne M$ and $K\ne N$ (if not, one of the exceptional cases (3a) or (3b) holds), $H$ and $K$ are Sylow $2$- or $3$-subgroups of $M$ and $N$, respectively. Now Lemma~\ref{WieHartSimple} and $H,K\in\sm_\pi(S)$ imply that $H,K$ are Sylow subgroups in $S$ and they are conjugate by the Sylow theorem.

Consider the case where $M$ and $N$ are non-isomorphic maximal subgroups of $S$. In this case, $|H|=|K|$ divides $(|M|,|N|)$ and, by the arguments similar above, we can consider that the numbers $|H|=|K|$ and $(|M|,|N|)$  are not  powers of primes.  Consider all possibilities for~$S$.

Let $S=L_2(q)$ where $q=2^p$, $p$~is a prime. Since $(|M|,|N|)$  is not a power of prime, by the information in Table~\ref{MaxL_2(2)}   we can consider that $M=E_q:C_{(q-1)}$  and
  $N=D_{2(q-1)}$. If  $2\in\pi$, then the $\pi$-Hall subgroups of $M$ and $N$ are non-isomorphic. Hence $2\notin\pi$. In this case, both $H$ and $K$ are abelian $\pi$-Hall subgroups of $S$ and they are conjugate by Lemma~\ref{NilpHall}.

 Let $S=L_2(q)$, where  $q=3^p$ and  $p$~is an odd prime. Note that the order of $E_q:C_{\frac{1}{2}(q-1)}$ is odd and the orders of $D_{q-1}$ and $D_{q+1}$ are not divisible by~3. Now it easy to see from Table~\ref{MaxL_2(3)} that the condition that $(|M|,|N|)$  is not a power of prime implies that we can consider that $M=E_q:C_{\frac{1}{2}(q-1)}$,
  $N=D_{q-1}$, and $2\notin\pi$. In this case, both $H$ and $K$ are abelian $\pi$-Hall subgroups of $S$ and they are conjugate by Lemma~\ref{NilpHall}.

 Let $S=L_2(q)$ where $q$~is a prime such that $q>3$ and $q^2+1\equiv 0\pmod 5$. The case when one of
$M$ and $N$ is $E_q:C_{\frac{1}{2}(q-1)}$ and other one  is $D_{q-1}$ can be argued similarly as the previous case, and as $(|M|,|N|)$  is not a power of prime, we can consider that $M\in\{A_4,S_4\}$. Hence $2,3\in\pi$. If $M=A_4$, then $H=M$ (if not, $|H|$ is a power of 2 or~3), but the maximal subgroups of $S$ of the other types  contain no subgroups isomorphic to~$A_4$. Hence $M=S_4$. Since the other maximal subgroups of $S$ do not contain subgroups isomorphic to $S_4$ and $H$ is not a power of a prime, we have that $H\cong S_3$. Since  a maximal subgroup $N$ of $S$ contains a subgroup $K$ isomorphic to $H\cong S_3$, it is easy to see from   Table~\ref{MaxL_2(q)} that $N=D_{q-\varepsilon}$ for some $\varepsilon\in\{+,-\}$ and in view of $K\in\Hall_\pi(N)$ we have that $|S|_3=|N|_3=3$. Hence $N$ is the normalizer of a Sylow 3-subgroup $P$ of $S$ and $K\in\Hall_\pi(N_S(P))$. But $H\cong K$ means that $H$ is also  a $\pi$-Hall subgroup of some (solvable) normalizer of a Sylow 3-subgroup $Q$ of $S$. Hence by 
the Sylow and Hall theorems, $H$ and $K$ are conjugate.

  Let $S=Sz(q)$, where  $q=2^p$, $p$~is an odd prime. This case can be investigated without essential changes as the case $S=L_2(q)$, $q=2^p$.

 Finally, let $S=L_3(3)$.  Since $(|M|,|N|)$  is not a power of prime and in view of the information in Table~\ref{MaxL_3(3)}, we can consider that $M=E_{3^2}:GL_2(3)$,
  $N=S_4$ and $\pi=\{2,3\}$. But $H\in\Hall_\pi(M)$ and $K\in\Hall_\pi(N)$ implies that $H=M\not\cong N=K$; a contradiction.
\medskip

Statement (4) is a straightforward consequence of~(3).
\qed\medskip

\subsection{
The pronormality of the $\pi$-submaximal subgroups\\ of minimal nonsolvable groups}\label{3part}

In order to complete the proof of Theorem~\ref{Th1}, we need to establish the pronormality of the $\pi$-maximal subgroups in minimal simple groups, that is, we need to show that every subgroup $H$ appearing in one of Tables~\ref{L_2(2),2_notin_pi}--\ref{L_3(3),2,3} is pronormal in the corresponding group $S\in\cal T$.

\begin{Pro}\label{Pronormality} {\it Let $S$ be a minimal simple group  and
$H\in\sm_\pi(S)$. Then $H$ is pronormal in $S$. }
\end{Pro}
\noindent{\sc Proof}. Let $g\in S$. We need to show that $H$ and $H^g$ are conjugate in $\langle H,H^g\rangle$. It is trivial if $\langle H,H^g\rangle=S$. Hence we can consider that $\langle H,H^g\rangle\leq M$ for some $M\lessdot S$. If $H\in\Hall_\pi(M)$, then $H,H^g\in\Hall_\pi(\langle H,H^g\rangle)$ and the solvability of $M$ implies that $H$ and  $H^g$ are conjugate in $\langle H,H^g\rangle$.

Thus, in view of the statement (2) of Proposition~\ref{ProofOfSubmaximality}, we only need to consider the cases (2a)--(2d) in this statement.
   \begin{itemize}
    \item[{\rm (2a)}] $S=L_2(q)$ where  $q$~is a prime such that $q>3$ and $q^2+1\equiv 0\pmod 5$, $2,3\in \pi$, $H=D_{(q-\varepsilon)_\pi}$ $\varepsilon\in\{+,-\}$,
       $\text{}\pi\cap\pi(q-\varepsilon)=\{2\}$, $q\equiv-\varepsilon7\pmod{48}$, and
          $H\cong D_8$ is a Sylow $2$-subgroup of $S$.
          \end{itemize}

          In this case, $H$ is pronormal as a Sylow subgroup of $S$.

          \begin{itemize}
    \item[{\rm (2b)}] $S=L_2(q)$ where  $q$~is a prime such that $q>3$ and $q^2+1\equiv 0\pmod 5$, $2,3\in \pi$, $H=D_{(q-\varepsilon)_\pi}
    $ $\varepsilon\in\{+,-\}$,
     $\pi\cap\pi(q-\varepsilon)=\{2,3\}$, $q\equiv\varepsilon7,\varepsilon31\pmod{72}$, and
     $H\cong S_3\cong D_6$ is a $\{2,3\}$-Hall subgroup of the normalizer of a Sylow $3$-subgroup of $S$.
         \end{itemize}

         By the Sylow theorem, one can consider that $H$ and $H^g$ are $\{2,3\}$-Hall subgroups of the same (solvable) normalizer of a  Sylow $3$-subgroup of~$S$. Then by the Hall theorem, $H$ and $H^g$ are conjugate in $\langle H,H^g\rangle$.

          \begin{itemize}
        \item[{\rm (2c)}] $S=L_3(3)$, $\pi\cap\pi(S)=\{2,3\}$, and $H=3_+^{1+2}:C_2^2$ is the normalizer of a Sylow $3$-subgroup of $S$. \end{itemize}

            In this case, $H$ and $H^g$ are conjugate in $K=\langle H,H^g\rangle$ since both $H$ and $H^g$ are normalizers of Sylow 3-subgroups of $K$.

          \begin{itemize}
    \item[{\rm (2d)}] $S=L_3(3)$, $\pi\cap\pi(S)=\{2,3\}$, and $H=GL_2(3)$.
  \end{itemize}
  It is easy see  that in this case $H$ contains a Sylow 2-subgroup $P$ of $S$. One can take some $x\in\langle H,H^g\rangle$ such that $P^x=P^g$. Then $gx^{-1}\in N_S(P)$. But $N_S(P)=P$ by \cite[Corollary]{Kond} and $g\in Px\subseteq \langle H,H^g\rangle$. Hence $H$ and $H^g$ are conjugate in $K=\langle H,H^g\rangle$.

\qed\medskip

\noindent{\sc Proof} of Corollary~\ref{Cor}. The corollary  is a straightforward consequence of Propositions~\ref{SubFactor}(4) and~\ref{Pronormality}.

%
%
%
%
%
%
%
%
%
%
%
%
%
%
%
%
%
%
%
%
%
%
%
%
%
%
%
%
%
%
%
%


\smallskip

\vspace{3cm}


\begin{thebibliography}{100}

\bibitem{Bray}
{J.~N.~Bray, D.~F.~Holt, C.~M.~Roney-Dougal}, The Maximal Subgroups of the Low-Dimensional Finite Classical Groups. Cambridge: Cambridge Univ. Press, 2013. 438~p.



\bibitem{Atlas}
{J.~H.~Conway, R.~T.~Curtis, S.~P.~Norton, R.~A.~Parker, R.~A.~Wilson,}
Atlas of Finite Groups. Oxford: Clarendon Press, 1985. 252~p.


\bibitem{Chun} S. A. Chunikhin,
\"{U}ber aufl\"{o}sbare Gruppen,  Mitt. Forsch.-Inst. Math. Mech. Univ. Tomsk, 2 (1938),  222--223.

\bibitem{DH}
{K.\,Doerk, T.\,Hawks,}
{Finite Soluble Groups},
Berlin, New York, Walter de Gruyter, 1992.

\bibitem{GuoBook}
{W.\,Guo,}
{The Theory of Classes of Groups},
Beijing, New York, Kluwer Acad. Publ., 2006.

\bibitem{GuoBook1}
{W.\,Guo,}
{Structure Theory of Canonical Classes of Finite Groups}, Berlin, Springer, 2015.


\bibitem{GR} W.\,Guo, D.\,O.\,Revin,   On the class of groups with pronormal Hall $\pi$-subgroups, Siberian Math. J., 55:3 (2014), 415--427.

\bibitem{GR1} W.\,Guo, D.\,O.\,Revin,   On maximal and submaximal ${\mathfrak X}$-subgroups, Algebra and logic, submitted.

\bibitem{GRV} W.\,Guo, D.\,O.\,Revin, E.\,P.\,Vdovin,  Confirmation for Wielandt's conjecture, J. Algebra, 434 (2015), 193--206.

\bibitem{Hart}
{ B.\,Hartley, }
A theorem of Sylow type for a finite groups,
{ Math.\ Z.}, { 122:4} (1971),  223--226.


\bibitem{Hart1}
{ B.\,Hartley, }
Helmut Wielandt on the $\pi$-structure of finite groups,
Helmut Wielandt: Mathematical Works, Vol. 1,
Group theory (ed. B. Huppert and H. Schneider, de Gruyter, Berlin, 1994), 511--516.


 \bibitem{Hall1}
{P. Hall,}
A note on soluble groups,
{ J. London Math. Soc.},  { 3} (1928), 98--105.

\bibitem{Hall2}
{P. Hall,}
A characteristic property of soluble groups,
{ J.\ London\ Math.\ Soc.}, {12} (1937), 198--200.



\bibitem{Hall} P. Hall,  Theorems like Sylow`s,
 Proc. London Math. Soc., 6:22 (1956), 286--304.




\bibitem{Hupp} B. Huppert, Endliche Gruppen. Berlin: Springer-Verlag,  1967.











\bibitem{Kleid}
{P.~B.~Kleidman,}
 A proof of the Kegel--Wielandt conjecture on subnormal subgroups., Ann. of Math. (2) 133 (1991), N2, 369--428.

 \bibitem{KL}
{P.~B.~Kleidman,  M.~Liebeck,} The Subgroup Structure of the Finite Classical Groups. Cambridge: Cambridge University Press, 1990. 303~p.


 \bibitem{Kond}
 A.~S.~Kondratiev, Normalizers of the Sylow 2-subgroups in finite simple groups, Math. Notes, 78:3 (2005), 338--346.




\bibitem{R4}  D.\,O. Revin, The $D_\pi$-property in finite simple groups, Algebra and Logic, 47:3 (2008), 210--227.



 \bibitem{DpiCl} D. O. Revin, The $D_\pi$-Property in a Class of Finite Groups, Algebra and Logic, 41:3 (2002), 187--206.





 \bibitem{Hall3'} D.\,O.\,Revin, E. P. Vdovin, Hall subgroups of finite groups, Contemporary Mathematics, 402 (2006), 229--265.




 \bibitem{NumbCl} D.\,O.\,Revin, E.\,P.\,Vdovin,  On the number of classes of conjugate Hall subgroups in finite simple groups, J. Algebra, 324:12 (2010),  3614--3652.



\bibitem{FrattArg} D.\,O.\,Revin, E.\,P.\,Vdovin, Frattini argument for Hall subgroups, J. Algebra, 414 (2014), 95--104.






 \bibitem{Shem} L. A. Shemetkov, Generalizations of Sylow's theorem, Siberian Math. J., 44:6 (2003), 1127--1132.

\bibitem{SuzCl}
 M.~Suzuki, On a class of doubly transitive groups, Ann. Math., 75:1 (1962), 105--145.



\bibitem{Suz} M.\,Suzuki, { Group Theory} II,  Springer-Verlag, New
York--Berlin--Heidelberg--Tokyo, 1986.



\bibitem{Thomp}
{J.~G.~Thompson},
 Nonsolvable finite groups all of whose local subgroups are solvable,  Bull. Amer. Math. Soc., 74 (1968), 383--437.




\bibitem{ConjCrit} E.\,P.\,Vdovin, D.\,O.\,Revin, A conjugacy criterion for Hall subgroups in finite groups, Siberian Math. J., 51:3 (2010), 402--409.


  \bibitem{ExCrit} D.\,O.\,Revin, E.\,P.\,Vdovin, An existence criterion for Hall subgroups of finite groups, J. Group Theory, 14:1 (2011), 93--101.

 \bibitem{Surv} E.\,P.\,Vdovin, D.\, O.\,Revin, Theorems of Sylow type, Russian Math. Surveys, 66:5 (2011), 829--870.






\bibitem{ProSimple}   E.\,P.\,Vdovin, D.\,O.\,Revin, Pronormality of Hall subgroups in finite simple groups, Siberian Math. J., 53:3 (2012), 419--430.


\bibitem{ProHall}   E.\,P.\,Vdovin, D.\,O.\,Revin, On the pronormality of Hall subgroups, Siberian Math. J., 54:1 (2013), 22--28.

\bibitem{ExProHall}   E.\,P.\,Vdovin, D.\,O.\,Revin, The existence of pronormal $\pi$-Hall subgroups in $E_\pi$-groups, Siberian Math. J., 56:3 (2015), 379--383.





\bibitem{Wie0}
{ H. Wielandt,}
Eine Verallgemeinerung der invarianten Untergruppen,
{ Math.\ Z.},  { 45} (1939), 209--244.

\bibitem{Wie}
{ H. Wielandt},
Zum Satz von Sylow,
{ Math.\ Z.},  { 60} (1954), N4. 407--408.


\bibitem{Wie2}
{ H. Wielandt},
Sur la Stucture des groupes compos\'es,
{ S\'eminare Dubriel-Pisot(Alg\`ebre et Th\'eorie des Nombres),} 17e
an\'ee, 10 pp. 1963/64. N17.

\bibitem{Wie3}
{ H. Wielandt},
Zusammenghesetzte Gruppen: H\"older Programm heute,
{ The Santa Cruz conf. on finite groups, Santa Cruz, 1979}. Proc.
Sympos. Pure
Math., { 37}, Providence RI: Amer. Math. Soc., 1980, 161--173.

\bibitem{Wie4}
{ H. Wielandt}
Zusammenghesetzte Gruppen endlicher Ordnung,
{ Vorlesung an der Universit\"at T\"ubingen im Wintersemester 1963/64}.
Helmut Wielandt: Mathematical Works, Vol. 1,
Group theory (ed. B. Huppert and H. Schneider, de Gruyter, Berlin, 1994), 607--516.

\bigskip
\end{thebibliography}
\end{document}